\documentclass[a4paper,10pt]{amsart}
\usepackage[utf8]{inputenc}
 \usepackage{setspace}
\usepackage[cmtip,arrow]{xy}
\usepackage{lscape}
\usepackage{latexsym}
\usepackage[activeacute,english]{babel}
\usepackage{graphicx}
\usepackage{amsmath}
\usepackage{amsfonts}
\usepackage{hyperref}
\hypersetup{colorlinks=true,linkcolor=blue,citecolor=red,linktocpage=true}

\usepackage{pst-node}

\usepackage{amsthm}
\usepackage{amssymb}
\usepackage{amscd}
\usepackage{mathrsfs}
\usepackage{pb-diagram,pb-xy}
\xyoption{all}
\usepackage{graphics}
\usepackage{stmaryrd}
\usepackage{yfonts,bbm}
\usepackage{tikz-cd}

\newtheorem{theorem}{Theorem}[section]
\newtheorem{lemma}[theorem]{Lemma}
\newtheorem{proposition}[theorem]{Proposition}
\newtheorem{corollary}[theorem]{Corollary}
\newtheorem{definition}[theorem]{Definition}
\newtheorem{example}[theorem]{Example}
\newtheorem{note}{Note}

\newtheorem{Remark}[theorem]{Remark}

\usepackage[square,sort,comma,numbers]{natbib}
\usepackage{hyperref}
\hypersetup{colorlinks=true,linkcolor=blue,citecolor=red,linktocpage=true}
\usepackage{url}
\usepackage{relsize}

\newcommand{\A}{\mathcal{A}}
\newcommand{\B}{\mathcal{B}}
\newcommand{\C}{\mathcal{C}}
\newcommand{\D}{\mathcal{D}}
\newcommand{\T}{\mathcal{T}}
\newcommand{\U}{\mathcal{U}}

\newcommand{\Z}{\mathbb{Z}}
\newcommand{\Ho}{\mathrm{Hom}}

\newcommand{\Dg}{\mathrm{DgMod}}

\title{Differential graded triangular matrix categories}
\author[M. L. S. Sandoval, V. Santiago, E. O. Velasco]{Martha Lizbeth Shaid Sandoval-Miranda\\ 
Valente Santiago-Vargas\\ Edgar Omar Velasco-P\'aez}

\thanks{The authors are grateful to the project PAPIIT-Universidad Nacional Aut\'onoma de M\'exico IN100520}
\dedicatory{}
\begin{document}
\maketitle

\begin{abstract}
This paper focuses on defining an analog of differential-graded triangular matrix algebra in the context of differential-graded categories. Given two  dg-categories $\U$ and $\T$ and $M \in \mathrm{DgMod}(\U \otimes \T^{op})$ we construct the differential graded triangular matrix category $\Lambda:=\left[\begin{smallmatrix}
\mathcal{T} & 0 \\ 
M & \mathcal{U}
\end{smallmatrix}\right]$. Our main result, is that there is  an equivalence of dg-categories  between the dg-comma category $\Big(\mathrm{DgMod}(\mathcal{T} ), \mathbb{G}\big(\mathrm{DgMod}(\mathcal{U})\big)\Big)$  and the category $\mathrm{DgMod}(\left[\begin{smallmatrix}
\mathcal{T} & 0 \\ 
M & \mathcal{U}
\end{smallmatrix}\right])$. This result is an extension of a well known result for artin algebras (see for example \cite[III.2]{AusBook}).
\end{abstract}

\section{Introduction}\label{sec:1}

The differential graded categories (dg-categories) and their dg-modules have played a fundamental role in mathematics for a long time. In the '70s and '80s of the last century,  dg-categories were the primary tool to study the matrix problems related to the representation theory of algebras. In addition, the development of their study has led to Drozd's proof of the tame-wild dichotomy (see, for example: \cite{Yuriy1}, \cite{Yuriy2}, \cite{Mark}). \\

In recent decades, the primary relevance of dg-categories has stemmed from a fundamental result due to B. Keller,  which asserts that each compactly generated algebraic triangulated  category is equivalent to the derived category of certain small dg-category, see \cite[Theorem 4.3]{Keller}.
This importance  increased even more when G. Tabuada  proved in  \cite[Theorem 0.1]{Tabuada} that the category of small dg-categories admits a model structure in which the weak equivalences are the quasi-equivalences; and so also, with the work of  B. T\"oen  (see \cite{Toen}), who studied in depth the associated homotopy category Ho(Dg-cat), showing in particular that it has an internal Hom and deriving several applications of this fact to homotopy theory and  algebraic geometry.\\

On the other hand, rings of the form $\left[\begin{smallmatrix}
T & 0 \\ 
M & U
\end{smallmatrix}\right]$ where $T$ and $U$ are rings and $M$ is a $T$-$U$-bimodule has often appeared in the study of the representation theory of artin rings and algebras (see for example \cite{Auslander1},\cite{Ringel1}, \cite{Gordon1}, \cite{Green1}, \cite{Green2}). Such rings, called {\it{triangular matrix rings}}, appear naturally in the study of homomorphic images of hereditary artin algebras and in the study of the decomposition of algebras and direct sum of two rings. Triangular matrix rings and their homological properties have been 
the subject of research for decades, see for example, \cite{Chase1}, \cite{Haghany1}, \cite{Michelena}, \cite{Skow1}, \cite{Ying1}.\\

It is well known that the category $\mathrm{Mod}\big( \left[\begin{smallmatrix}
T & 0 \\ 
M & U
\end{smallmatrix}\right] \big)$ of all left modules over the triangular matrix ring $\left[\begin{smallmatrix}
T & 0 \\ 
M & U
\end{smallmatrix}\right]$ is equivalent to the comma category $(M{\otimes_{T}}- \downarrow 1_{\mathrm{Mod}(U)}),$ whose
objects are triples $(A,f, B)$ where $f : M \otimes_{T} A \rightarrow B $ is a morphism of $U$-modules,
and the morphisms are pairs $(\alpha, \beta)$ with $\alpha :A \rightarrow A'$ and $\beta: B \rightarrow B'$ such that the following diagram commutes
$$\xymatrix{ M \otimes_{T} A \ar[r]^{1_{M} \otimes \alpha} \ar[d]_{f} & M \otimes_{T} A' \ar[d]^{f'}\\
B \ar[r]_{\beta}& B'}.$$ For more details, see for example, \cite[III.2]{AusBook}.\\

Note that in the context of  differential graded algebras, triangular matrix rings are also studied. For instance, D. Maycock's work on dg-algebra  $\left[\begin{smallmatrix}
T & 0 \\ 
M & U
\end{smallmatrix}\right],$ where $T$ and $U$ are dg-algebras. In \cite{Maycock}, D. Maycock constructed a recollement between the derived categories of certain dg-algebras, and he studied derived equivalences between the upper triangular matrix  dg-algebras.\\

Given a ring $R$, there exists a preadditive category  associated to $R$, the category  $\mathcal{R}$ that has only one object, say $\ast$, and 
$\mathrm{Hom}_\mathcal{R}(\ast,\ast)=R$.  We can consider  the abelian category of all additive covariant functors $F:\mathcal{R}\longrightarrow \textbf{Ab}$ denoted by $\mathrm{Mod}(\mathcal{R})$, where $\textbf{Ab}$ denotes the category of abelian groups. The evaluation  functor
$\mathrm{ev}_{\ast}: \mathrm{Mod}(\mathcal{R})\rightarrow \mathrm{Mod}(R)$, $\mathrm{ev}_{\ast}(F)=F(\ast)$ induces an equivalence of categories. So, given a preadditive category $\mathcal{C},$ we can think of $\mathcal{C}$ as a kind of ring with several objects, and we can consider  its category of $\mathcal{C}$-modules, $\mathrm{Mod}(\mathcal{C}),$  which consists of all the additive covariant functors  $F:\mathcal{C}\longrightarrow \textbf{Ab}$. Following this ideas, in \cite{LeOS},
A. Le\'on-Galeana, M. Ortiz-Morales, and V. Santiago, 
defined the triangular matrix categories $\Lambda=\left[\begin{smallmatrix} \mathcal{T} & 0 \\ 
M & \mathcal{U}  \end{smallmatrix}\right],$ where $\mathcal{T}$ and $\mathcal{U}$ are preadditive categories, and $M:\mathcal{U}\otimes_{\mathbb{Z}}\mathcal{T}^{op}\longrightarrow \bf{Ab}$ is a bifunctor. Several properties of the category $\mathrm{Mod}(\Lambda)$ were studied there. \\

In this paper we define the analogous of the  differential graded triangular matrix algebra in the context of dg-categories. Given two dg-categories $\mathcal{U}$ and $\mathcal{T}$ and $M \in \mathrm{DgMod}(\U \otimes \T^{op})$ a dg-bifunctor, we construct the triangular matrix category $\Lambda := \left[\begin{smallmatrix}
\T & 0 \\ 
M & \U
\end{smallmatrix}\right]$; and we prove that $\left[\begin{smallmatrix}
\mathcal{T} & 0 \\ 
M & \mathcal{U}
\end{smallmatrix}\right]$ is a dg-category, see Proposition \ref{matrizdgcat}.\\ 
Then, we prove that there is  an equivalence of dg-categories  between the dg-comma category $\Big(\mathrm{DgMod}(\mathcal{T} ), \mathbb{G}\big(\mathrm{DgMod}(\mathcal{U})\big)\Big),$  and the dg-category $\mathrm{DgMod}(\left[\begin{smallmatrix}
\mathcal{T} & 0 \\ 
M & \mathcal{U}
\end{smallmatrix}\right]).$ This is, Theorem \ref{equivalencia Lambda y coma}:
\begin{center}
\it The differential graded functor
$\mathfrak{F}: \Big( \Dg(\T), \big(\mathbb{G}\Dg(\U)\big) \Big) \longrightarrow \Dg(\Lambda) $\\
is a dg-equivalence of dg-categories.
\end{center}

It is worth of noticing that this description of $\mathrm{Mod}(\Lambda)$ as a comma category is useful because it is easier to study the properties of $\mathrm{Mod}(\Lambda)$ through the categories $\mathrm{Mod}(\mathcal{U})$ and $\mathrm{Mod}(\mathcal{T})$ involved in this comma category (see for example \cite[Proposition 2.3 ]{AusBook}) in p. 75).\\

\section{Preliminaries}
In this section we present some basic notations, mainly taken from \cite{Keller} and \cite{Keller2} (see also \cite{Yeku1} for a textbook), which will be used all throughout these paper. Throughout this section $K$ will be an arbitrary commutative ring. Recall that a differential graded (dg) $K$-module is a graded $K$-module $M= \bigoplus_{n \in \Z} M^{n}$, together with a graded $K$-linear morphism $d_{M}: M\longrightarrow M$ of degree $+1$ such that $d_{M}\circ d_{M} = 0$. The category which will be  denoted by $\mathrm{DgMod}(K)$ has as objects the dg $K$-modules; and
given $M$ and $N$ two dg $K$-modules the set of morphisms from $M$ to $N$ in the category $\mathrm{DgMod}(K)$ is by definition the following graded $K$-module: 
$$\Ho_{\mathrm{DgMod}(K)}(M, N ):= \bigoplus_{n \in \Z} \Ho_{\mathrm{DgMod}(K)}^{n}(M, N ),$$
where $ \Ho_{\mathrm{DgMod}(K)}^{n}(M, N )$ consists of the graded $K$-linear morphisms $f:M \longrightarrow N$ of degree $n$, i.e., such that $f(M^{i})\subseteq N^{i+n}$ for all $i \in \Z$.
Moreover, each space of morphisms $\Ho_{\mathrm{DgMod}(K)}(M, N) $ has a structure of differential graded  $K$-module given by the  differential $d: \Ho_{\mathrm{DgMod}(K)}(M, N) \longrightarrow \Ho_{\mathrm{DgMod}(K)}(M, N)$, which is a graded $K$-linear morphism of degree $1$ such that $d \circ d = 0$ and is defined by the rule:
\begin{equation}\label{Homdgestructure}
d(\alpha) = d_{N} \circ \alpha -(-1)^{|\alpha|} \alpha \circ d_{M}
\end{equation}
where $|?|$ denotes the degree, whenever $\alpha$ is a homogeneous element
of the set $\Ho_{\mathrm{DgMod}(K)}(M, N)$.\\
Now,  if $M$ and $N$ are dg $K$-modules, the tensor product $M \otimes N := M \otimes_{K} N$ also becomes an object of $\mathrm{DgMod}(K)$ where the grading is given by $(M \otimes N )^{n} = \bigoplus_{i+j=n} M^{i}\otimes_{K} N^{j}$ and the differential $d_{M \otimes N} : M \otimes N \longrightarrow M \otimes N$ by the rule
\begin{equation} \label{tendgestructure}
d_{M \otimes N} (m \otimes n) = d_{M}(m) \otimes n + (-1)^{|m|}m \otimes d_{N}(n),
\end{equation}
for all homogeneous elements $m \in M$ and $n \in N$. All throughout this paper, we use the unadorned symbol $\otimes$ to denote $\otimes_{K}$.\\

\begin{definition}
 A differential graded category or a $\textbf{dg-category}$ is a $K$-category $\C $ such that:
 \begin{itemize}
 \item[(a)] $\Ho_{\mathcal{C}}(X,Y)$ $\in \mathrm{DgMod}(K)$   for all $ X,Y \in \mathcal{C}$.
 \item[(b)] The composition function  
 $$\theta_{X,Y,Z}: \Ho_{\mathcal{C}}(Y,Z) \otimes_{K} \Ho_{\mathcal{C}}(X,Y) \longrightarrow \Ho_{\mathcal{C}}(X,Z)$$
is a morphism of degree zero in $\mathrm{DgMod}(K)$  which commutes with the differentials.
 \end{itemize}
 \end{definition}
 
\begin{example}
\begin{enumerate}
\item[(a)] A dg-category with only one object can be identified with a dg-algebra, i.e., a graded $K$-algebra with a differential $d$ such that it satisfies Leibniz's rule $d(fg) = d(f)g + (-1)^{|f|}fd(g)$
for f, g homogeneous elements. Conversely, each dg-algebra $B$ can be viewed as a dg-category with only one object.
\item[(b)] The category $\mathrm{DgMod}(K)$ is the main example of a dg-category.
\end{enumerate}
\end{example}

If $\A$ is a dg-category, then the $\textbf{opposite dg-category}$ $\A^{op}$ has the same class of objects as $\A$ and the differential on morphisms $d:\A^{op}(A, B)= \A(B, A)\longrightarrow \A(B, A) = \A^{op}(A, B)$ is the same as in $\A$, that is, $d_{\A^{op}(A, B)}(f^{op}):=(d_{\A(B,A)}(f))^{op}$. The composition of morphisms in $\A^{op}$ is given as $\beta^{op} \circ  \alpha^{op} = (-1)^{|\alpha||\beta|}(\alpha \circ \beta)^{op}$, where we use the superscript $op$ to emphasize that a morphism is viewed in $\A^{op}$.\\
Analogously to the tensor product of $K$-modules we recall the definition of the tensor product of two differential graded categories.
\begin{definition}
Let $\A,$ $\B$ be dg-categories. The $\textbf{tensor product of dg-categories}$, which we denote by $\A \otimes \B$ is defined as follows:
 \begin{itemize}
  \item[(a)] Obj($\mathcal{A} \otimes \B$):= $Obj(\A) \times Obj(\B)$.
  \item[(b)]$\Ho_{\A \otimes \B}((X,Y),(X',Y')):= \Ho_{\A}(X,X') \otimes \Ho_{\B}(Y, Y')$ for all  $(X,Y),$ $ (X',Y') \in Obj(\A \otimes \B)$.
\end{itemize}
The composition 
  $\circ:\Ho_{\A \otimes \B}((X',Y'),(X'',Y'')) \times \Ho_{\A \otimes \B}((X,Y),$ $(X',Y')) \longrightarrow \Ho_{\A \otimes \B}((X,Y),(X'',Y''))$ is given as:
  $$g \circ f= (\alpha_{2} \otimes \beta_{2}) \circ (\alpha_{1} \otimes \beta_{1})= (-1)^{|\beta_{2}| |\alpha_{1}|} (\alpha_{2}\alpha_{1}) \otimes (\beta_{2}\beta_{1}),$$ 
 for each $f=\alpha_{1} \otimes \beta_{1}$ and $g=\alpha_{2} \otimes \beta_{2}$ with
$\alpha_{1}\in \Ho_{\A}(X,X')$, $\alpha_{2}\in \Ho_{\A}(X',X'')$, $\beta_{1} \in \Ho_{\B}(Y,Y')$ and  $\beta_{2} \in \Ho_{\B} (Y',Y'')$ homogenous elements.
\end{definition}

\begin{Remark}\label{diferentialtensorcat}
It is important to note that the existence of the tensor product of two dg-categories is essentially due to the tensor product of differential graded $K$-modules given in Equation \ref{tendgestructure}.\\
Indeed, for a pair of objects $(X,Y)$, $(X',Y')\in \mathcal{A}\otimes \mathcal{B}$ we get that the set of morphisms $\mathrm{Hom}_{\mathcal{A}\otimes \mathcal{B}}((X,Y),(X',Y'))$ is a dg K-module with differential
$$\delta:\mathrm{Hom}_{\mathcal{A}\otimes \mathcal{B}}((X,Y),(X',Y'))\longrightarrow  \mathrm{Hom}_{\mathcal{A}\otimes \mathcal{B}}((X,Y),(X',Y'))$$
defined as follows: for $\alpha\otimes\beta\in \mathrm{Hom}_{\mathcal{A}\otimes \mathcal{B}}((X,Y),(X',Y'))=\Ho_{\A}(X,X') \otimes \Ho_{\B}(Y, Y')$  a homogeneous element we set
$$\delta(\alpha\otimes \beta)=d_{\mathcal{A}(X,X')}(\alpha)\otimes \beta+(-1)^{|\alpha|}\alpha\otimes d_{\mathcal{B}(Y,Y')}(\beta).$$
\end{Remark}

\begin{definition}
Let  $\A$, $\B$ be two  dg-categories.  A  \textit{graded functor} $F: \A \longrightarrow \B$   is a  $K$-linear functor such that $F(\Ho_{\A}^{n} (A, A'))\subseteq \Ho_{\B}^{n} (F(A), F(A'))$  for  each  $n \in \mathbb{Z}$ and  $A, A' \in Obj(\A)$.
 \end{definition}
 
 \begin{definition}
 Let $\A$, $\B$ two dg-categories. A  $\textbf{differential graded functor}$ or $\textbf{dg-functor}$ is a  graded functor  $F: \A \longrightarrow \B$  which commutes with  the differential, i.e, $F(d_{\A}(\alpha))= d_{\B}(F(\alpha))$ for all homogeneous morphism $\alpha$ in $\A$.
 \end{definition}
 One of the most important examples of functor requires the following lemma.
  \begin{lemma} \label{AxB dg funtor} 
Let $\A$, $\B$  and $\C$ be  dg-categories and  let $F: \A \otimes \B \longrightarrow \C$ be  an assignment on objects $(A,B) \rightsquigarrow F(A,B)$ and an assignment  on homogeneous morphism  $\alpha \otimes \beta \rightsquigarrow F(\alpha \otimes \beta)$ such that $|F(\alpha \otimes \beta)|= |\alpha|+ |\beta|$.Then the following conditions are equivalent.
\begin{itemize}
\item[(a)]The given mapping defines a dg-functor $F: \mathcal{A} \otimes \mathcal{B} \longrightarrow \mathcal{C}$.
\item[(b)]The following conditions are satisfied:
\begin{itemize}
\item[(b1)]For each fixed object  $A \in \A$, the assignment $B \rightsquigarrow F(A,B)$ on objects and $\beta \rightsquigarrow F(1_{A} \otimes \beta)$ on morphisms, define a dg-functor $G: \mathcal{B} \longrightarrow \mathcal{C}$.
\item[(b2)] For each fixed object $B \in \B$, the assignment $A \rightsquigarrow F(A,B)$ on objects and  $\alpha \rightsquigarrow F(\alpha \otimes 1_{B})$ on morphisms define a dg-functor $H: \mathcal{A} \longrightarrow \mathcal{C}$.
\item[(b3)]For all homogeneous morphisms $\alpha: A \longrightarrow A'$ and  $\beta: B \longrightarrow B'$, in  $\A$ and $\B$, respectively, there is the equality
$$(-1)^{|\alpha| |\beta |} F(1_{A'} \otimes \beta) \circ F(\alpha \otimes 1_{B})= F(\alpha \otimes \beta)= F(\alpha \otimes 1_{B'})\circ F(1_{A} \otimes \beta).$$
\end{itemize}
\end{itemize}
\end{lemma}
\begin{proof}
See \cite[Lemma 1.1]{Saorin}.
\end{proof}

\begin{example}\label{Hom(A,-) dg funtor}
\begin{enumerate}
\item [(a)] 
Let $\A$ be an  dg-category and  $\mathrm{DgMod}(K)$  the category  of differential graded $K$-modules. Given $A \in \A$, then 
$$\Ho_{\A}(A,-): \A \longrightarrow \mathrm{DgMod}(K)$$ is a dg-functor, where $\Ho_{A}(A,f)\!:\!\! \Ho_{\A}(A,A') \rightarrow \Ho_{\A} (A,B')$ is given by  $\Ho_{\A}(A,f)(j)= f \circ j$ for $f$, $j$  homogeneous elements.

\item [(b)] Let $\A^{op}$ be a dg-category, and let $\mathrm{DgMod}(K)$ be the category  of differential graded $K$-modules. If $A \in \A^{op}$, then $$\Ho_{\A^{op}}(-,A): \A^{op} \longrightarrow \mathrm{DgMod}(K)$$  is a dg-functor, where  $\Ho_{\A^{op}}(f^{op},A): \Ho_{\A^{op}}(B',A) \longrightarrow \Ho_{\A^{op}} (A',A)$ is given by $\Ho_{\mathcal{A}^{op}}(f,A)(j)= (-1)^{|f||j|} j \circ f$ for each  pair  $f$, $j$ of homogeneous elements.
\end{enumerate}
\end{example}

Using the Proposition \ref{AxB dg funtor} and examples  \ref{Hom(A,-) dg funtor} it follows that $\Ho_{\mathcal{A}}(-,-): \A^{op} \otimes \A  \longrightarrow \mathrm{DgMod}(K)$ is a dg-functor given for each  $(A,A') \in Obj(\A^{op} \otimes \A)$  as  $\Ho_{\mathcal{A}}(-,-)(A,A'):= \Ho_{\mathcal{A}}(A,A')$ and for   $\alpha: A \longrightarrow B$ and $ \alpha': A' \longrightarrow B'$  homogeneous morphism  in $\A$  we have that $\Ho_{\A}(\alpha^{op} \otimes \alpha')$ is defined  as follows:  $\Ho_{\A}(\alpha^{op} \otimes \alpha')(f):= (-1)^{|\alpha| (|\alpha'|+|f|)} \alpha' \circ f \circ \alpha$ for each homogeneous element $f \in \Ho_{\A}(B,A')$. 

Let us remember that the $\textbf{category of small dg-categories}$, which we denote by $\mathbf{dg}$-$\mathbf{Cat_{K}}$, is the category that has as objects all the small dg-categories, and as morphisms the dg-functors between them. To avoid sizing problems, it is requested that  $\mathrm{dg}$-$\mathrm{Cat}_{K}$ consists only of small dg-categories, fixed to a given universe, similar to $\mathrm{Cat}$ being the category of all small categories. To be able to work with dg-categories we require the notion of natural transformation in the dg context.

Now, we consider $\C$ and $\D$ two small dg-categories. In order to show that  the category of covariant dg-functors from $\mathcal{C}$ to $\mathcal{D}$ possesses a natural structure of differentially graded category, we need the following definition.

 \begin{definition}\label{dgnaturaltrans}
Let $\mathcal{C}$ and $\mathcal{D}$ be  small dg-categories, let $F,G: \mathcal{C} \longrightarrow \mathcal{D}$ be two  dg-functors. A \textbf{dg-natural transformation of degree n}  denoted by  $\eta: F \longrightarrow G$,  is a family of morphisms $\{ \eta_{X}: F(X) \longrightarrow G(X)\}_{X \in \mathcal{C}}$ such that $\eta_{X} \in \Ho_{\mathcal{D}} ^{n} (F(X), G(X))$ for each $X \in \mathcal{C}$ and for all $f \in \Ho_{\mathcal{C}} ^{m} (X, Y)$ homogeneous morphism of degree $m$,  the following diagram  commutes up to the sign $(-1)^{nm}$ 
$$\xymatrix{ F(X) \ar[r]^{\eta_{X}} \ar[d]_{F(f)}& G(X) \ar[d]^{G(f)} \\
  F(Y) \ar[r]_{\eta_{Y}} & G(Y), }$$
i.e.,  $G(f)\circ \eta_{X} = (-1)^{nm} \eta_{Y}\circ F(f)$. Let us denote by $\mathrm{DgNat}^{n}(F,G)$ the set of all the dg-natural transformation of degree $n$ and we define:
$$\mathrm{DgNat}(F,G):=\bigoplus_{n\in \mathbb{Z}}\mathrm{DgNat}^{n}(F,G),$$
as the set of all the $\textbf{dg}$-$\textbf{natural transformations}$ from $F$ to $G$.
\end{definition} 

The sign $(-1)^{mn}$ appearing in the above definition is due to the fact that a dg-category is an enriched category over the category of complexes of $K$-modules and the latter is a symmetric monoidal category, whose symmetry $c: \mathrm{DgMod}(K) \otimes \mathrm{DgMod}(K) \longrightarrow \mathrm{DgMod}(K) \otimes \mathrm{DgMod}(K)$ is given in homogenous elements by $c( x \otimes y )= (-1)^{|x||y|}y \otimes x$ for more information see \cite{Kelly} section 1.4.\\
Let $\mathcal{C}$ and $\mathcal{D}$ be  small dg-categories. Taking into account the dg-natural transformations one can now construct a new category whose objects will be dg-functors going from $\C$ to $\D$ and whose morphisms are the natural dg-transformations $\mathrm{DgNat}(F,G)$. In this way,  the $\textbf{category of dg-covariant functors}$ denoted by $\mathrm{DgFun}(\C,\D)$ is defined as follows. 
\begin{enumerate}
\item [(a)] The class of objects are the dg-covariant functors $F: \C \longrightarrow \D$.
\item [(b)] For $F,G  \in \mathrm{DgFun}(\C,\D)$ the set of morphisms will be by definition the set of natural dg-transformations from $F$ to $G$, i.e,
$$\Ho_{\mathrm{DgFun}(\C,\D)}(F,G):= \mathrm{DgNat}(F,G).$$

\end{enumerate}
Moreover each space of dg-natural transformations $\mathrm{DgNat}(F,G)$ has a structure of differential graded $K$-module given by the differential $d:\mathrm{DgNat}(F,G)\longrightarrow \mathrm{DgNat}(F,G)$ which is a graded $K$-linear morphism of degree $1$ such that $d \circ d = 0$ and is defined by the rule:

\begin{equation}\label{diferentialNat}
d(\eta)_{X}:=d_{\mathrm{Hom}_{\mathcal{D}}(F(X),G(X))}(\eta_{X}),\,\,\, \forall \eta\in \mathrm{DgNat}^{n}(F,G)
\end{equation}
and for all $X\in \mathcal{C}$, where the morphism $d_{\mathrm{Hom}_{\mathcal{D}}(F(X),G(X))}:\mathrm{Hom}_{\mathcal{D}}(F(X),G(X))\longrightarrow \mathrm{Hom}_{\mathcal{D}}(F(X),G(X))$ is the differential in the set $\mathrm{Hom}_{\mathcal{D}}(F(X),G(X))$ coming from the fact that $\mathcal{D}$ is a dg-category.

The following proposition is easy to verify.

\begin{proposition}\label{Hom(C,D) dg cat}
Let $\mathcal{C}$, $\mathcal{D}$ be  small dg-categories. Then the category $\mathrm{DgFun}(\C,\D)$  form a dg-category.
 \end{proposition}

If we take $K$ as the dg-category with the single object $\{\star\}$ as the unit
object, we have that the category $\mathrm{dg}$-$\mathrm{Cat}_{K}$ is a symmetric tensor category, i.e. we have the
adjunction
$$ \Ho_{\mathrm{dg}-\mathrm{Cat}_{K}}(\A \otimes \B, \C)\cong \Ho_{\mathrm{dg}-\mathrm{Cat}_{K}}(\A,\mathrm{Hom}_{\mathrm{dg}-\mathrm{Cat}_{K}}(\B,\C))$$
for $\A,\B,\C$ dg-categories. The pair $(\otimes,\Ho)$ makes $\mathrm{dg}$-$\mathrm{Cat}_{K}$ into a closed symmetric monoidal category (see \cite{Kelly}). This structure will be important for the remainder of this work  since the immediate consequence of the above fact is that there exists an isomorphism of dg-categories

\begin{equation}\label{bifuntor2variables}
\mathrm{DgFun}(\A \otimes \B, \C)\cong \mathrm{DgFun}(\A,\mathrm{DgFun}(\B,\C)).
\end{equation}
Note that this  isomorphism is related to the Lemma \ref{AxB dg funtor}.\\
Just as dg-algebras generalize to dg-categories, the same game can be played with dg-modules and this is crucial for the rest of the paper. We will introduce dg-modules.

\begin{definition}
Let $\C$ be a small dg-category. We  define a  {\textbf{left dg $\C$-module}} to be a dg-functor  $M: \C \longrightarrow \mathrm{DgMod}(K)$ while a right dg $\C$-module is a dg-functor  $N: \C^{op} \longrightarrow \mathrm{DgMod}(K)$.
\end{definition}

Now, we define the dg-category of left dg $\mathcal{C}$-modules.
\begin{definition}
Let $\C$ be a small dg-category. The category of left dg $\C$-modules denoted by $\mathrm{DgMod}(\C)$ has all dg $\C$-modules as objects and morphisms of dg-functors as the dg-natural transformations. That is,
$$\mathrm{DgMod}(\C):=\mathrm{DgFun}(\mathcal{C},\mathrm{DgMod}(K)).$$
\end{definition}

By Proposition \ref{Hom(C,D) dg cat}, we have that $\mathrm{DgMod}(\C)$ is a dg-category such that for any dg-functors $F,G: \C \longrightarrow \mathrm{DgMod}(K)$

\begin{equation}\label{morfismosdgmod}
\Ho_{\mathrm{DgMod}(\C)}(F,G)= \bigoplus_{n \in \Z} \Ho^{n}_{\mathrm{DgMod}(\C)}(F,G),
\end{equation}

where by definition we set $\Ho^{n}_{\mathrm{DgMod}(\C)}(F,G):=\mathrm{DgNat}^{n}(F,G)$ (see Definition \ref{dgnaturaltrans}). Moreover,  $\Ho_{\mathrm{DgMod}(\C)}(F,G)$ has a structure of differential graded $K$-module given by the differential $d: \Ho_{\mathrm{DgMod}(\C)}(F,G)\longrightarrow  \Ho_{\mathrm{DgMod}(\C)}(F,G)$ which is a graded $K$-linear morphism of degree $1$ such that $d \circ d = 0$ and is defined by the rule:

\begin{equation}\label{dgNatmod}
d(\eta)_{X}:=d_{\mathrm{Hom}_{\mathrm{DgMod}(K)}(F(X),G(X))}(\eta_{X}),\,\,\, \forall \eta\in  \Ho^{n}_{\mathrm{DgMod}(\C)}(F,G)
\end{equation}
and for all $X\in \mathcal{C}$, where the morphism 
$$d_{\mathrm{Hom}_{\mathrm{DgMod}(K)}(F(X),G(X))}\!:\!\mathrm{Hom}_{\mathrm{DgMod}(K)}(F(X),G(X))\rightarrow \mathrm{Hom}_{\mathrm{DgMod}(K)}(F(X),G(X))$$ is the differential  defined in the Equation  \ref{Homdgestructure}.\\

Now we recall the Yoneda emdedding in the context of differential graded categories, for more details see Section 4.3 in  \cite{Takeda} or Equation (2.34) in p. 34 in  \cite{Kelly} for the general context of enriched categories.

\begin{definition}(Dg-Yoneda's embedding)\label{dgYoneda}
Let $\mathcal{C}$ be a dg-category. Then the dg-Yoneda embedding
$$Y:\mathcal{C}\longrightarrow \mathrm{DgMod}(\mathcal{C}^{op})$$
given as $Y(C):=\mathrm{Hom}_{\mathcal{C}}(-,C)$  is a dg-funtor.
\end{definition}

\section{Differential Graded Triangular Matrix Categories}

Before we begin to construct the category of matrices, we construct the so-called {\it comma category in the differential graded context.}

\begin{definition}\label{commadgcat}
Let $\mathcal{C}$ and $\mathcal{D}$ be two dg-categories, and $F:\mathcal{C}\longrightarrow \mathcal{D}$ be a dg-functor. We define the $\textbf{comma category}$ $\big(\mathcal{D},F(\mathcal{C})\big)$ as follows:
\begin{enumerate}
\item [(a)] The objects of $\big(\mathcal{D},F(\mathcal{C})\big)$ are morphisms $f:D\longrightarrow F(C)$ in $\mathcal{D}$ with $f\in \mathrm{Hom}_{\mathcal{D}}^{0}(D,F(C))$ and $d_{\mathcal{D}(D,F(C))}(f)=0$. The object
$f:D\longrightarrow F(C)$ will be denoted by $(D,f,C)$.

\item [(b)] Given two objects $f:D\longrightarrow F(C)$ and $f:D'\longrightarrow F(C')$ in $\big(\mathcal{D},F(\mathcal{C})\big)$ a morphism from $f$ to $f'$ is a pair of morphisms $\alpha:D\longrightarrow D'$ and $\beta:C\longrightarrow C'$ such that the following diagram commutes in $\mathcal{D}$:
$$\xymatrix{D\ar[r]^{\alpha}\ar[d]_{f} & D'\ar[d]^{f'}\\
F(C)\ar[r]^{F(\beta)} & F(C').}$$
\end{enumerate}
The composition is given as follows: let $(\alpha,\beta):(D,f,C)\longrightarrow (D',f',C')$ and $(\alpha',\beta'):(D',f',C')\longrightarrow (D'',f'',C'')$ morphisms in  $\big(\mathcal{D},F(\mathcal{C})\big)$ we set:
$$(\alpha',\beta')\circ (\alpha,\beta):=(\alpha'\circ \alpha,\beta'\circ\beta).$$
\end{definition}

\begin{proposition}\label{dgcommaprop}
Let $\mathcal{C}$ and $\mathcal{D}$ be two dg-categories, and $F:\mathcal{C}\longrightarrow \mathcal{D}$ be a dg-functor. Then the comma category $\big(\mathcal{D},F(\mathcal{C})\big)$ is a dg-category.
\end{proposition}
\begin{proof}
It is easy to see that  $\big(\mathcal{D},F(\mathcal{C})\big)$ is a category.\\
Now, let us see that $\mathrm{Hom}_{\big(\mathcal{D},F(\mathcal{C})\big)}\Big((D,f,C),(D',f',C')\Big)$ 
is a dg $K$-module. For simplicity we will use the notation $\mathrm{Hom}_{\big(\mathcal{D},F(\mathcal{C})\big)}(f,f')$ to denote the previous set of morphisms.\\
We define
$$\mathrm{Hom}_{\big(\mathcal{D},F(\mathcal{C})\big)}^{n}\big(f,f'\big):=\left \{(\alpha,\beta)\mid \alpha\in \mathcal{D}^{n}(D,D'),\,\, \beta\in \mathcal{C}^{n}(C,C'),\,\, f'\circ \alpha=F(\beta)\circ f\right\}.$$
Then we have that $\mathrm{Hom}_{\big(\mathcal{D},F(\mathcal{C})\big)}\big(f,f'\big)=\bigoplus_{n\in \mathbb{Z}}\mathrm{Hom}_{\big(\mathcal{D},F(\mathcal{C})\big)}^{n}\big(f,f'\big).$\\
Now, we define $\delta:\mathrm{Hom}_{\big(\mathcal{D},F(\mathcal{C})\big)}\big(f,f'\big)\longrightarrow \mathrm{Hom}_{\big(\mathcal{D},F(\mathcal{C})\big)}\big(f,f'\big)$
as follows: for $(\alpha,\beta)\in \mathrm{Hom}_{\big(\mathcal{D},F(\mathcal{C})\big)}^{n}\big(f,f'\big)$ we set
$$\delta(\alpha,\beta):=(d_{\mathcal{D}(D,D')}(\alpha),d_{\mathcal{C}(C,C')}(\beta)).$$
Using the fact that $F$ is a dg-functor and the formula given in Equation \ref{Homdgestructure}, we can see that the following diagram commutes
$$\xymatrix{D\ar[rr]^{d_{\mathcal{D}(D,D')}(\alpha)}\ar[d]_{f} & & D'\ar[d]^{f'}\\
F(C)\ar[rr]^{F(d_{\mathcal{C}(C,C')}(\beta))} & & F(C').}$$
This proves that $\delta(\alpha,\beta)=(d_{\mathcal{D}(D,D')}(\alpha),d_{\mathcal{C}(C,C')}(\beta))\in \mathrm{Hom}_{\big(\mathcal{D},F(\mathcal{C})\big)}^{n+1}\big(f,f'\big)$. It is straightforward that $\delta\delta=0$ and hence $\mathrm{Hom}_{\big(\mathcal{D},F(\mathcal{C})\big)}(f,f')$ is a dg K-module with differential $\delta$.

Now, let us see that the composition is a morphism of degree zero in $\mathrm{DgMod}(K)$, that is, that the following diagram commutes
$$(\ast):\xymatrix{\mathrm{Hom}(f',f'')\otimes\mathrm{Hom}(f,f')\ar[rr]^{\Theta}\ar[d]_{\Delta} & & \mathrm{Hom}(f,f'')\ar[d]^{d_{\mathrm{Hom}(f',f'')}}\\
\mathrm{Hom}(f',f'')\otimes\mathrm{Hom}(f,f')\ar[rr]^{\Theta} & & \mathrm{Hom}(f,f'')}$$
Indeed, let  $(\alpha,\beta)\in \mathrm{Hom}(f,f')$ and $(\alpha',\beta')\in \mathrm{Hom}(f',f'')$ be homogenous elements with $|(\alpha,\beta)|=|\alpha|=|\beta|$ and $|(\alpha',\beta')|=|\alpha'|=|\beta'|$.
On one hand,  we have that

\begin{align*}
& (\Theta\circ \Delta)\Big((\alpha',\beta')\otimes (\alpha,\beta)\Big)=\\
&= \Theta\Big(d_{\mathrm{Hom}(f',f'')}(\alpha',\beta')\otimes (\alpha,\beta)+(-1)^{|(\alpha',\beta')|}(\alpha',\beta')\otimes d_{\mathrm{Hom}(f,f')}(\alpha,\beta) \Big)\\
&= d_{\mathrm{Hom}(f',f'')}(\alpha',\beta')\circ (\alpha,\beta)+(-1)^{|(\alpha',\beta')|}(\alpha',\beta')\circ d_{\mathrm{Hom}(f,f')}(\alpha,\beta) \Big)\\
& = \Big(d_{\mathcal{D}(D,D'')}(\alpha'\circ \alpha),\,\,\, d_{\mathcal{C}(C,C'')}(\beta'\circ \beta)\Big).
\end{align*}
On the other hand,
\begin{align*}
(d_{\mathrm{Hom}(f,f'')}\circ \Theta)\Big((\alpha',\beta')\otimes (\alpha,\beta)\Big) & =d_{\mathrm{Hom}(f,f'')}\Big((\alpha'\circ \alpha,\beta'\circ \beta)\Big)\\
& = \Big(d_{\mathcal{D}(D,D'')}(\alpha'\circ \alpha),\,\,\, d_{\mathcal{C}(C,C'')}(\beta'\circ \beta)\Big).
\end{align*}
This proves that the diagram $(\ast)$ commutes and hence the comma category $\big(\mathcal{D},F(\mathcal{C})\big)$ is a dg-category.
\end{proof}

 \begin{proposition}\label{dg funtores de  M}
Let $\U$ and $\T$ be dg-categories and $M \in \mathrm{DgMod}(\U \otimes \T^{op})$.  Then, there  exist  two  covariant dg-functors 
\begin{center}
$E': \T^{op} \longrightarrow \mathrm{DgMod}(\U)$\\
$Q: \U \longrightarrow \mathrm{DgMod}(\T^{op}).$
\end{center}
Moreover, $E'$ induces a dg-functor $E: \T \longrightarrow \mathrm{DgMod}(\U)^{op}$.
\end{proposition}
\begin{proof}
By Equation \ref{bifuntor2variables}, we have that existence of the dg-functors $E'$ and $Q$. Since, we will use the descriptions of the functors $E$ and $Q$, we recall the definition of such functors.\\
$(a)$. For $T \in \T$, we defined  a covariant  functor  $E'(T):= M_{T}: \U \longrightarrow \mathrm{DgMod}(K)$ as follows:
\begin{enumerate}
\item[(a1)] $M_{T}(U):=M(U,T)$,  for all $U \in \U$.
\item[(a2)] $M_{T}(u):= M(u \otimes 1_{T}^{op})$,  for all $u \in \Ho_{\U}(U, U')$.
\end{enumerate}
By Lemma \ref{AxB dg funtor}, we have that $E'(T)$ is a dg-funtor.\\
Now, given  a homogeneous morphism  $t:T \longrightarrow T'$ in  $\T$  we set  $E'(t^{op}):= \overline{t}: M_{T'} \longrightarrow M_{T}$ where $\overline{t}:= \{ [\overline{t}]_{U}: M_{T'}(U) \longrightarrow M_{T}(U)\}_{U \in \U}$ with  $[\overline{t}]_{U \in \U}:=M(1_{U} \otimes t^{op}): M(U,T') \longrightarrow M(U,T)$.  This defines a covariant dg-functor $E':\mathcal{T}^{op}\longrightarrow \mathrm{DgMod}(\U)$. Then we can define a dg-functor $E:\T \longrightarrow \mathrm{DgMod}(\U)^{op}$ as $E(t):=(E'(t^{op}))^{op}$.\\
 $(b)$. We defined  $Q: \U \longrightarrow \mathrm{DgMod}(\T^{op})$ first in objects. Let $U \in \U$ be, then  $Q(U):= M_{U}: \T^{op} \longrightarrow \mathrm{DgMod}(k)$ is  a functor  defined as follow   
\begin{enumerate}
\item[(b1)] $M_{U}(T):=M(U,T)$  for all $T \in \T^{op}$.
\item[(b2)] $M_{U}(t):= M(1_{U} \otimes t^{op})$ for each $t \in \Ho_{\T^{op}}(T, T')$.
\end{enumerate}
By Lemma \ref{AxB dg funtor}, $Q(U)$ is a dg-functor. Now for a homogeneous morphism $u:U \longrightarrow U'$ en $\U$ we define a dg-natural transformation $Q(u):= \overline{u}: M_{U} \longrightarrow M_{U'}$ where $\overline{u}= \{ [\overline{u}]_{T}: M_{U}(T) \longrightarrow M_{U'}(T)\}_{T \in \T}$ and $[\overline{u}]_{T \in \T}:=M(u \otimes 1_{T}^{op}): M(U,T) \longrightarrow M(U',T)$. 
\end{proof}

Now, we recall that $\mathrm{DgMod}(\U)$ is a dg-category that consists of the dg-functors
$M:\U\longrightarrow \mathrm{DgMod}(K)$ and given $M$ and $N$ in  $\mathrm{DgMod}(\U)$  we have that the space of morphisms  $\mathrm{Hom}_{\mathrm{DgMod}(\U)}(M, N)$ consists of the dg-natural transformations (see Equation \ref{morfismosdgmod}).

\begin{definition}\label{definiG}
We define a covariant dg-functor $\mathbb{G}: \mathrm{DgMod}(\U) \longrightarrow \mathrm{DgMod}(\T)$ as follows.  Let $Y: \mathrm{DgMod}(\U) \longrightarrow \mathrm{DgMod}(\mathrm{DgMod}(\U)^{op})$ be the Yoneda dg-functor  $Y(B):= \Ho_{\mathrm{DgMod}(\U)}(-,B)$ (see Definition \ref{dgYoneda})  and consider the dg-functor $I:\mathrm{DgMod}(\mathrm{DgMod}(\U)^{op}) \longrightarrow \mathrm{DgMod}(\T)$, induced by  $E: \T \longrightarrow \mathrm{DgMod}(\U)^{op}$ of the  Proposition \ref{dg funtores de  M}.  We set  $\mathbb{G}:= I \circ Y: \mathrm{DgMod}(\U) \longrightarrow \mathrm{DgMod}(\T)$.
\end{definition}
\begin{note} \label{Obs de G}
In detail, we have that the following holds.
 \begin{enumerate}
 \item For  $B \in \mathrm{DgMod}(\U)$, $\mathbb{G}(B)(T)=\Ho_{\mathrm{DgMod}(\U)}(M_{T},B)$ for all $T \in \T$. Now, for  $t: T \longrightarrow T'$ a homogeneous morphism in $\T$ we have $\mathbb{G}(B)(t):\Ho_{\mathrm{DgMod}(\U)}(M_{T},B) \longrightarrow \Ho_{\mathrm{DgMod}(\U)}(M_{T'},B)$  and  for each $\eta: M_{T} \longrightarrow B$  homogeneous morphism of degree $|\eta|$  we have $\mathbb{G}(B)(t)(\eta):= \Ho (\overline{t}, B)(\eta)$ $= (-1)^{|\eta||\overline{t}|} \eta \circ \overline{t}$, where $\overline{t}:= E'(t^{op}).$
 
 \item  Let $\varepsilon: B \longrightarrow B'$ be a  dg-natural transformation of dg $\U$-modules of degree $|\varepsilon|$. We assert that  $\mathbb{G}(\varepsilon):\mathbb{G}(B) \longrightarrow \mathbb{G}(B')$ is a dg-natural transformation in $\mathrm{DgMod}(\T)$ of degree $|\varepsilon|$. Indeed, for each  $T \in \T$ we have that $[\mathbb{G}(\varepsilon)]_{T}:= \Ho_{\mathrm{DgMod}(\U)}(M_{T}, \varepsilon):  \Ho_{\mathrm{DgMod}(\U)}(M_{T}, B)\longrightarrow  \Ho_{\mathrm{DgMod}(\U)}(M_{T}, B')$. Also remember that the following  diagram  commutes up to the sign $(-1)^{|\varepsilon| |t|}$:
 $$\xymatrix{ \Ho_{\mathrm{DgMod}(\U)}(E(T),B)\ar[rr]_{\Ho_{\mathrm{DgMod}(\U)} (E(T),\varepsilon)} \ar[dd]_{\Ho_{\mathrm{DgMod}(\U)} (E(t),B)}& &\Ho_{\mathrm{DgMod}(\U)}(E(T),B') \ar[dd]_{\Ho_{\mathrm{DgMod}(\U)} (E(t),B')} \\
  & & \\
   \Ho_{\mathrm{DgMod}(\U)}(E(T'),B)\ar[rr]_{\Ho_{\mathrm{DgMod}(\U)}(E(T'),\varepsilon)}& & \Ho_{\mathrm{DgMod}(\U)}(E(T'),B').}$$
   for each homogeneous morphism  $t: T \longrightarrow T'$ in  $\T$. As a result $$\mathbb{G}(\varepsilon)=\big \{ [\mathbb{G}(\varepsilon)]_{T}:= \Ho_{\mathrm{DgMod}(\U)}(E(T),\varepsilon): \mathbb{G}(B)(T) \longrightarrow \mathbb{G}(B')(T) \big \}_{T \in \T}$$
is dg-natural transformation of degree $|\varepsilon|$.
 \end{enumerate}
\end{note}
 
\begin{Remark}
The previous note shows that
 $$\mathbb{G}\Big(\mathrm{Hom}_{\mathrm{DgMod}(\U)}^{n}(B,B')\Big)\subseteq \mathrm{Hom}_{\mathrm{DgMod}(\T)}^{n}(\mathbb{G}(B),\mathbb{G}(B')).$$
In fact, this must happen since $\mathbb{G}$ is a dg-functor according to Definition  \ref{definiG}.
\end{Remark}
 
Let $1_{\mathrm{DgMod}(\T)}:\mathrm{DgMod}(\T) \longrightarrow \mathrm{DgMod}(\T)$ and $\mathbb{G}:\mathrm{DgMod}(\U) \longrightarrow \mathrm{DgMod}(\T)$ be  dg-functors. Therefore,  from the Definition \ref{commadgcat},  we have  the dg-comma category $\Big( \Dg(\T), \mathbb{G}\big(\Dg(\U)\big) \Big)$ whose objects are  the triples  $(A,f,B)$ with $A \in \mathrm{DgMod}(\T)$, $B \in \mathrm{DgMod}(\U)$ and $f:A \longrightarrow \mathbb{G}(B)$ a morphism of differential graded  $\T$-modules of degree zero and such that $d_{\mathrm{Hom}(A,\mathbb{G}(B))}(f)=0$. That is, $f_{T}\in \mathrm{Hom}_{\mathrm{DgMod}(K)}^{0}\big(A(T),\mathbb{G}(B)(T)\big)$ for each $T\in \mathcal{T}$ and 
for each homogeneous morphism $t \in \Ho_{\T}^{n}(T,T')$ of degree $n$,  the following diagram  commute  in $\mathrm{DgMod}(K)$:
\begin{equation}\label{digramacoma}
\xymatrix{ A(T) \ar[rr]^{f_{T}} \ar[dd]^{A(t)}& & \mathbb{G}(B)(T) \ar[dd]^{\mathbb{G}(B)(t)=\mathrm{Hom}(\overline{t},B)} \\
& & \\
A(T') \ar[rr]_{f_{T'}}& &  \mathbb{G}(B)(T').}
\end{equation}

By Definition \ref{commadgcat}, a morphism of degree $k$  between two objects  $(A,f,B)$ and  $(A',f',B')$ is a pair  of  homogenous morphisms  $(\alpha, \beta)$ with  $\alpha: A \longrightarrow A'$ a morphism of dg $\T$-modules and  $\beta: B \longrightarrow  B$ a morphism  of dg $\U$-modules  such $k=|\alpha|=|\beta|$ and such that  the following  diagram commutes 
 $$\xymatrix{ A \ar[r]^{\alpha} \ar[d]_{f}& A'  \ar[d]^{f'} \\
  \mathbb{G}(B) \ar[r]_{\mathbb{G}(\beta)}&  \mathbb{G}(B').}$$
Note  that, since  $f: A \longrightarrow \mathbb{G}(B)$ is a morphism differential graded $\T$-modules of degree zero,  for each homogeneous morphism $t \in \Ho_{\T}^{n}(T,T')$ of degree $n$,  the diagram  in Equation \ref{digramacoma}, commutes  in $\mathrm{DgMod}(K).$ That is, $\mathbb{G}(B)(t) \circ f_{T}= f_{T'}\circ A(t)$. Then for  $x \in A(T) $ and homogenous element of degree $r=|x|$ we have $f_{T}(x) \in (\mathbb{G}(B)(T))^{r}:=\Ho_{\mathrm{DgMod}(\U)}^{r}(M_{T}, B)$. Then,
$$f_{T}(x)= \{ [f_{T}(x)]_{U}: M_{T}(U) \longrightarrow B(U) \}_{U \in \U},$$ 
where $[f_{T}(x)]_{U}\in \mathrm{Hom}_{\mathrm{DgMod}(K)}^{r}\big(M_{T}(U),B(U)\big)$ for each $U\in \mathcal{U}$.

\begin{lemma}\label{productos dg matrices}
Let $M \in \mathrm{DgMod}(\mathcal{U} \otimes \mathcal{T}^{op})$ be and $f:A \longrightarrow \mathbb{G}(B)$  a morphism in $\mathrm{DgMod}(\T)$ of degree zero and such that  $d_{\mathrm{Hom}(A,\mathbb{G}(B))}(f)=0$.  Let $U \in \U$ and $T \in \T$, for $m \in M(U,T)$ and  $x\in A(T)$  homogeneous elements  we set: 
$$m \cdot x:= (-1)^{|x||m|}[f_{T}(x)]_{U}(m) \in B(U)^{|x|+|m|}$$
(this  product can be  defined for each morphism of degree zero $f: A \longrightarrow \mathbb{G}(B)$). We extend this product linearly as follows: if  $x= \sum_{n}x_{n}\in A(T)$  and $m= \sum_{k}m_{k}\in M(U,T)$ are sum of homogeneous elements,  we define
$$m \cdot x:= \sum_{n,k \in \mathbb{Z}}m_{k} \cdot x_{n}.$$
 \begin{itemize}
\item[(a)] For  $m \in M(U,T')$, $t \in \Ho_{\mathcal{T}}(T,T')$ and  $x \in A(T)$ homogenous elements  we set:
$$m \bullet t:= (-1)^{|t||m|} M(1_{U} \otimes t^{op})(m), \hspace{1cm} t \ast x := A(t)(x).$$
We extend this definition by linearlity as follows: if  $x= \sum_{n}x_{n}\in A(T)$, $m= \sum_{k}m_{k}\in M(U,T)$  and $t=\sum_{i}t_{i}\in \Ho_{\mathcal{T}}(T,T')$ are a sum of homogeneous elements we define: $m \bullet t:= \sum_{i,k \in \mathbb{Z}}m_{k} \bullet t_{i}$ and $t \ast x := A(t)(x).$
Then 
$$(m \bullet t)\cdot x=m \bullet (t \ast x),$$
for every $x\in A(T)$, $m\in M(U,T)$  and $t\in \Ho_{\mathcal{T}}(T,T')$.

\item[(b)] For  $m \in M(U,T)$, $u \in \Ho_{\mathcal{U}}(U,U')$ and $z \in B(U)$ elements (not necessarely homogenous) we set
$$u \bullet m:=M(u \otimes 1_{T})(m), \hspace{2cm} u \diamond z := B(u)(z).$$
Then  we have that $(u \bullet m) \cdot x= u \diamond (m \cdot x)$.
\item[(c)] Let $m_{1} \in M(U',T')$, $m_{2} \in M(U,T)$, $t \in \Ho_{\mathcal{T}}(T,T')$ and  $u \in \Ho_{\mathcal{U}}(U,U')$  elements (not necessarely homogenous). For  $x \in A(T)$  (not necessarely homogenous)  we have that
$$(m_{1} \bullet t + u \bullet m_{2}) \cdot x=( m_{1} \bullet t) \cdot x + (u \bullet m_{2}) \cdot x= m_{1} \cdot(t \ast x)+  u \diamond (m_{2} \cdot x).$$
\end{itemize}
\end{lemma}
\begin{proof}
By the definition of $m \cdot x$. It is enough to verify the conditions in homogeneous elements.
\begin{itemize}
\item[(a)] Since $f: A\longrightarrow \mathbb{G}(B)$ is a morphism of dg $\mathcal{T}$-module of degree zero, for each  $t \in \Ho_{\mathcal{T}}^{n}(T,T')$  the  diagram in Equation \ref{digramacoma}, commutes. That is, for $x \in A(T)$ we have that
\begin{align*}
f_{T'}\left( A(t)(x) \right) = \left(  \mathbb{G}(B)(t)\circ f_{T}\right)(x) & = \left( \Ho_{\mathrm{DgMod}(\mathscr{U})}(\overline{t},B) \circ f_{T}\right)(x)\\
& =(-1)^{|f_{T}(x)||\overline{t}|} f_{T}(x)\circ \overline{t}\\
& = (-1)^{|x||t|} f_{T}(x)\circ \overline{t}.
\end{align*}
Hence, for  $U \in \U$ we have that
$$  \big[ f_{T'} \left( A(t)(x) \right) \big]_{U}=(-1)^{|x||t|}[f_{T}(x)]_{U}\circ [\overline{t}]_{U}= (-1)^{|x||t|}[f_{T}(x)]_{U} \circ M(1_{U} \otimes t^{op}), $$
 and thus
$$ \big[ f_{T'} \left( A(t)(x) \right) \big]_{U} (m) =  (-1)^{|x||t|} [f_{T}(x)]_{U}( M(1_{U} \otimes t^{op})(m)).$$

Now, we have that
\begin{align*}
(m \bullet t)\cdot x & = \left(  (-1)^{|m||t|} M(1_{U}\otimes t^{op})(m)\right) \cdot x\\
& = (-1)^{|m||t|}(-1)^{|x|(|m|+|t|)}[f_{T}(x)]_{U}(M(1_{U} \otimes t^{op})(m))\\
& =  (-1)^{|m||t|}(-1)^{|x||m|}(-1)^{|x||t|}[f_{T}(x)]_{U}(M(1_{U} \otimes t^{op})(m))
\end{align*}
 where the the sign of the second last equality is due to the fact that M is a dg-funtor, so that 
$|M(1_{U} \otimes t^{op})|=|1_{U} \otimes t^{op}|=|t|$ and  in this way,  the morphism of dg k-modules $M(1_{U} \otimes t^{op})$ send the element $m$ to the component of degree $|m|+|t|$.\\
On the other hand,
\begin{align*}
m\cdot (t \ast x) & = m \cdot (A(t)(x))\\
& = (-1)^{|m|(|t|+|x|)}[f_{T}(A(t)(x))]_{U}(m) \\
& = (-1)^{|m||t|}(-1)^{|m||x|}[f_{T}(A(t)(x))]_{U}(m)\\
& =(-1)^{|m||t|}(-1)^{|m||x|} (-1)^{|x||t|} [f_{T}(x)]_{U}( M(1_{U} \otimes t^{op})(m)).
\end{align*}
where the second equality is because  $|A(t)(x)|=|t|+|x|$. Therefore, we conclude that $(m \bullet t)\cdot x=m \bullet (t \ast x)$.

\item[(b)] The proof of item $(b)$ and $(c)$ is similar to $(a)$.
\end{itemize}
\end{proof}
Now we define the differential graded triangular matrix category, although it is a definition, we must verify that it indeed form a category.
 \begin{definition}\label{Lamda definicion}
 We define the  {\bf{differential graded triangular matrix category }} $\Lambda=\left[\begin{smallmatrix} \mathcal{T} & 0 \\ 
M & \mathcal{U}  \end{smallmatrix}\right]$ as follows. 
\begin{itemize}
\item[(a)]  The objects of this dg-category are matrices   $ \left[\begin{smallmatrix} T & 0 \\ 
M & U  \end{smallmatrix}\right]$ with $T \in Obj(\mathcal{T})$ and  $U \in Obj(\mathcal{U})$.
\item[(b)] Given a pairs of objects  $ \left[\begin{smallmatrix} T & 0 \\ 
M & U  \end{smallmatrix}\right]$, $\left[\begin{smallmatrix} T' & 0 \\ 
M & U'  \end{smallmatrix}\right]$ in $\Lambda$  we define
$$\Ho_{\Lambda} \Big( \left[\begin{smallmatrix} T & 0 \\ 
M & U  \end{smallmatrix}\right], \left[\begin{smallmatrix} T' & 0 \\ 
M & U'  \end{smallmatrix}\right]  \Big) := \left[\begin{smallmatrix}
\Ho_{\mathcal{T}}(T,T') & 0\\
M(U',T) & \Ho_{\mathcal{U}}(U,U')
\end{smallmatrix}\right].$$
\end{itemize} 
The composition is given by
$$\circ : \left[\begin{smallmatrix}
\mathcal{T}(T',T'') & 0\\
M(U'',T') & \mathcal{U}(U',U'')
\end{smallmatrix}\right] \times  \left[\begin{smallmatrix}
\mathcal{T}(T,T') & 0\\
M(U',T) & \mathcal{U}(U,U')
\end{smallmatrix}\right] \longrightarrow  \left[\begin{smallmatrix}
\mathcal{T}(T,T'') & 0\\
M(U'',T) & \mathcal{U}(U',U'')
\end{smallmatrix}\right] $$
$$\Big( \left[\begin{smallmatrix}
t_{2} & 0 \\
m_{2} & u_{2}
\end{smallmatrix}\right],
\left[\begin{smallmatrix}
t_{1} & 0 \\
m_{1} & u_{1}
\end{smallmatrix}\right] \Big)
\longmapsto 
\left[\begin{smallmatrix}
t_{2}\circ t_{1} & 0 \\
m_{2} \bullet t_{1} + u_{2} \bullet m_{1} & u_{2} \circ u_{1}
\end{smallmatrix}\right].$$
We recall that for homogenous elements we have that $m_{2} \bullet t_{1}:= (-1)^{|m_{2}||t_{1}|}M(1_{U''} \otimes t_{1}^{op})(m_{2})$ and  $u_{2} \bullet m_{1}:= M(u_{2} \otimes 1_{T})(m_{1})$.
\end{definition}

\begin{lemma}\label{Lambda categoria}
The composition defined above is associative, and for an object
$\left[\begin{smallmatrix}
T & 0\\
M & U
\end{smallmatrix}\right]
 \in \Lambda$,  the identity morphism is given by 
$1_{\left[\begin{smallmatrix}
T & 0 \\ 
M & U
\end{smallmatrix}\right]} :=
\left[\begin{smallmatrix}
1_{T} & 0\\
0 & 1_{U}
\end{smallmatrix}\right]
$
\end{lemma}
\begin{proof}
It follows from Lemma \ref{productos dg matrices}.
\end{proof}

\begin{lemma}\label{Md=dM}
Let $M \in \mathrm{DgMod}(\U \otimes \T^{op})$ be. Then for  a homogeneous element $\alpha\otimes\beta^{op}\in \mathrm{Hom}_{\mathcal{U}\otimes \mathcal{T}^{op}}((U,T'),(U',T))$ we have that:

\begin{align*}
& M\Big(d_{\mathcal{U}(U,U')}(\alpha)\otimes \beta^{op}\Big)+(-1)^{|\alpha|}M\Big(\alpha\otimes d_{\mathcal{T}(T,T')}(\beta)^{op}\Big)=\\
&= d_{M(U',T)}\circ M(\alpha\otimes\beta^{op})-(-1)^{|\alpha|+|\beta|}M(\alpha\otimes\beta^{op})\circ d_{M(U,T')}.
\end{align*}
\end{lemma}
\begin{proof}
Since $M$ is a dg-functor, we get that $M$ commutes with the differential. That is,  for a homogenous element
$\alpha\otimes\beta^{op}\in \mathrm{Hom}_{\mathcal{U}\otimes \mathcal{T}^{op}}((U,T'),(U',T))$ we have that 
$$M(\delta(\alpha\otimes \beta^{op}))=d(M(\alpha\otimes \beta^{op})).$$
That is, we have the following commutative diagram for elements homogenous in the set  $\mathrm{Hom}_{\mathcal{U}\otimes \mathcal{T}^{op}}((U,T'),(U',T))$:
$$\xymatrix{\mathrm{Hom}_{\mathcal{U}\otimes \mathcal{T}^{op}}((U,T'),(U',T))\ar[rr]^{M}\ar[d]_{\delta} & & \mathrm{Hom}_{\mathrm{DgMod}(K)}\big(M(U,T'),M(U',T)\big)\ar[d]^{d}\\
\mathrm{Hom}_{\mathcal{U}\otimes \mathcal{T}^{op}}((U,T'),(U',T))\ar[rr]^{M} & & \mathrm{Hom}_{\mathrm{DgMod}(K)}\big(M(U,T'),M(U',T)\big).}$$

Using  the description of $\delta$ and $d$ given in Remark \ref{diferentialtensorcat} and Equation \ref{Homdgestructure}, we obtain that for a homogeneous element $\alpha\otimes\beta^{op}\in \mathrm{Hom}_{\mathcal{U}\otimes \mathcal{T}^{op}}((U,T'),(U',T))$ we have that

\begin{align*}
& M\Big(d_{\mathcal{U}(U,U')}(\alpha)\otimes \beta^{op}\Big)+(-1)^{|\alpha|}M\Big(\alpha\otimes d_{\mathcal{T}(T,T')}(\beta)^{op}\Big)=\\
&=M(\delta(\alpha\otimes \beta^{op})\\
& = d(M(\alpha\otimes \beta^{op}))\\
&= d_{M(U',T)}\circ M(\alpha\otimes\beta^{op})-(-1)^{|\alpha|+|\beta|}M(\alpha\otimes\beta^{op})\circ d_{M(U,T')},
\end{align*}
since $|M(\alpha\otimes \beta^{op})|=|\alpha|+|\beta|$.
\end{proof}

\begin{lemma}\label{similarreglaLeibniz}
Let
$$\left[
\begin{smallmatrix}
t_{1} & 0 \\
m_{1} & u_{1}
\end{smallmatrix}\right]\in \mathrm{ Hom}_{\mathbf{\Lambda}}\left (\left[ \begin{smallmatrix}
T & 0 \\
M & U
\end{smallmatrix} \right] ,  \left[ \begin{smallmatrix}
T' & 0 \\
M & U'
\end{smallmatrix} \right]  \right)  = \left[ \begin{smallmatrix}
\mathrm{Hom}_{\mathcal{T}}(T,T') & 0 \\
M(U',T) & \mathrm{Hom}_{\mathcal{U}}(U,U')
\end{smallmatrix} \right]$$

$$\left[
\begin{smallmatrix}
t_{2} & 0 \\
m_{2} & u_{2}
\end{smallmatrix}\right]\in \mathrm{ Hom}_{\mathbf{\Lambda}}\left (\left[ \begin{smallmatrix}
T' & 0 \\
M & U'
\end{smallmatrix} \right] ,  \left[ \begin{smallmatrix}
T'' & 0 \\
M & U''
\end{smallmatrix} \right]  \right)  = \left[ \begin{smallmatrix}
\mathrm{Hom}_{\mathcal{T}}(T',T'') & 0 \\
M(U'',T') & \mathrm{Hom}_{\mathcal{U}}(U',U'')
\end{smallmatrix} \right]$$

Then, the following equalities hold.
\begin{enumerate}
\item [(a)]  $d_{M(U'',T)}(m_{2}\bullet t_{1})=d_{M(U'',T')}(m_{2})\bullet t_{1}+(-1)^{|m_{2}|}m_{2}\bullet d_{\mathcal{T}(T,T')}(t_{1})$,

\item [(b)] $d_{M(U'',T)}(u_{2}\bullet m_{1}) = d_{\mathcal{T}(U',U'')}(u_{2})\bullet m_{1}+(-1)^{|u_{2}|}u_{2}\bullet d_{M(U',T)}(m_{1})$.
\end{enumerate}
\end{lemma}
\begin{proof}
$(a)$.  Consider $1_{U''}\otimes  t_{1}^{op}:(U'',T')\longrightarrow (U'',T)$ a morphism in $\mathcal{U}\otimes \mathcal{T}^{op}$. Since $d_{\mathcal{U}(U'',U'')}(1_{U''})=0$ and $|1_{U''}|=0$, by Lemma \ref{Md=dM}, we have that
$M\Big(1_{U''}\otimes d_{\mathcal{T}(T,T')}(t_{1})^{op}\Big)=d_{M(U'',T)}\circ M(1_{U''}\otimes  t_{1}^{op})-(-1)^{|t_{1}|}M(1_{U''}\otimes  t_{1}^{op})\circ d_{M(U'',T')},$
since $|M(1_{U''}\otimes  t_{1}^{op})|=|t_{1}|$.\\
Then for $m_{2}\in M(U'',T')$,  and $r:= (|t_{1}|+1)\cdot |m_{2}|$ we have that
\begin{align*}
& m_{2}\bullet d_{\mathcal{T}(T,T')}(t_{1})=(-1)^{(|t_{1}|+1)|m_{2}|}M\Big(1_{U''}\otimes d_{\mathcal{T}(T,T')}(t_{1})^{op}\Big)(m_{2})\\
& =(-1)^{r}M\Big(1_{U''}\otimes d_{\mathcal{T}(T,T')}(t_{1})^{op}\Big)(m_{2})\\
& =(-1)^{r}\Big(\Big(d_{M(U'',T)}\circ M(1_{U''}\otimes  t_{1}^{op})\Big)(m_{2})-(-1)^{|t_{1}|}\Big(M(1_{U''}\otimes  t_{1}^{op})\circ d_{M(U'',T')}\Big)(m_{2})\Big)\\
& =(-1)^{|m_{2}|}d_{M(U'',T)}(m_{2}\bullet t_{1})-(-1)^{r}(-1)^{|t_{1}|}\Big((-1)^{|t_{1}|(|m_{2}|+1)}d_{M(U'',T')}(m_{2})\bullet t_{1}\Big)\\
& =(-1)^{|m_{2}|}d_{M(U'',T)}(m_{2}\bullet t_{1})-\Big((-1)^{|t_{1}|(|m_{2}|+1)+|t_{1}|+(|t_{1}|+1)|m_{2}|}d_{M(U'',T')}(m_{2})\bullet t_{1}\Big)\\
& =(-1)^{|m_{2}|}d_{M(U'',T)}(m_{2}\bullet t_{1})-\Big((-1)^{|t_{1}|(|m_{2}|+1)+|t_{1}|+(|t_{1}|+1)|m_{2}|}d_{M(U'',T')}(m_{2})\bullet t_{1}\Big)\\
& =(-1)^{|m_{2}|}d_{M(U'',T)}(m_{2}\bullet t_{1})-\Big((-1)^{|t_{1}||m_{2}|+|t_{1}|+|t_{1}|+|t_{1}||m_{2}|+|m_{2}|}d_{M(U'',T')}(m_{2})\bullet t_{1}\Big)\\
& =(-1)^{|m_{2}|}d_{M(U'',T)}(m_{2}\bullet t_{1})-\Big((-1)^{|m_{2}|}d_{M(U'',T')}(m_{2})\bullet t_{1}\Big).
\end{align*}
Hence $d_{M(U'',T)}(m_{2}\bullet t_{1})=d_{M(U'',T')}(m_{2})\bullet t_{1}+(-1)^{|m_{2}|}m_{2}\bullet d_{\mathcal{T}(T,T')}(t_{1})$. This proves $(a)$.\\
$(b)$. The proof is similar to $(a)$.
\end{proof}

\begin{proposition}\label{matrizdgcat}
Let $\mathcal{T}$, $\mathcal{U}$ be two dg-categories and $M \in \mathrm{DgMod}(\mathcal{U} \otimes\mathcal{T}^{op})$. Consider $T,T' \in \mathcal{T}$ and $U, U' \in \mathcal{U}$.  Then
$$\Ho_{\Lambda} \Big( {\left[ \begin{smallmatrix} T & 0 \\ 
M & U  \end{smallmatrix}\right]},{\left[\begin{smallmatrix} T' & 0 \\ 
M & U'  \end{smallmatrix}\right]}  \Big) := \left[\begin{smallmatrix}
{\Ho_{\mathcal{T}}(T,T')} & 0\\
M(U',T) &{ \Ho_{\mathcal{U}}(U,U')}
\end{smallmatrix}\right]$$
has a differential graded K-module structure and the composition 
$$\circ : \left[\begin{smallmatrix}
\mathcal{T}(T',T'') & 0\\
M(U'',T') & \mathcal{U}(U',U'')
\end{smallmatrix}\right] \otimes  \left[\begin{smallmatrix}
\mathcal{T}(T,T') & 0\\
M(U',T) & \mathcal{U}(U,U')
\end{smallmatrix}\right] \longrightarrow  \left[\begin{smallmatrix}
\mathcal{T}(T,T'') & 0\\
M(U'',T) & \mathcal{U}(U',U'')
\end{smallmatrix}\right] $$
$$\big( \left[\begin{smallmatrix}
t_{2} & 0 \\
m_{2} & u_{2}
\end{smallmatrix}\right],
\left[\begin{smallmatrix}
t_{1} & 0 \\
m_{1} & u_{1}
\end{smallmatrix}\right] \big)
\longmapsto 
\left[\begin{smallmatrix}
t_{2}\circ t_{1} & 0 \\
m_{2} \bullet t_{1} + u_{2} \bullet m_{1} & u_{2} \circ u_{1}
\end{smallmatrix}\right]$$
is a morphism of dg K-modules of degree zero. In particular, $\Lambda=\left[\begin{smallmatrix} 
\mathcal{T} & 0\\
M & \mathcal{U}
\end{smallmatrix}\right]$ is a dg K-category.
\end{proposition}
\begin{proof}
Let $\mathcal{T}$ and $\mathcal{U}$ be two dg-categories and $M \in \mathrm{DgMod}(\U \otimes \T^{op})$, that is $M:\U \otimes \T^{op}\longrightarrow \mathrm{DgMod}(K)$ is a dg-functor.
Firstly, let us see that 
$$\Ho_{\Lambda} \Big( {\left[ \begin{smallmatrix} T & 0 \\ 
M & U  \end{smallmatrix}\right]},{\left[\begin{smallmatrix} T' & 0 \\ 
M & U'  \end{smallmatrix}\right]}  \Big) := \left[\begin{smallmatrix}
{\Ho_{\mathcal{T}}(T,T')} & 0\\
M(U',T) &{ \Ho_{\mathcal{U}}(U,U')}
\end{smallmatrix}\right]$$
is a dg $K$-module. As a $K$-module we have that
$$\Ho_{\Lambda} \big( \left[\begin{smallmatrix} T & 0 \\ 
M & U  \end{smallmatrix}\right], \left[\begin{smallmatrix} T' & 0 \\ 
M & U'  \end{smallmatrix}\right]  \big) = \bigoplus_{n \in \mathbb{Z}} \Ho_{\Lambda}^{n} \big( \left[\begin{smallmatrix} T & 0 \\ 
M & U  \end{smallmatrix}\right], \left[\begin{smallmatrix} T' & 0 \\ 
M & U'  \end{smallmatrix}\right]  \big),$$
with
$\Ho_{\Lambda}^{n} \big( \left[\begin{smallmatrix} T & 0 \\ 
M & U  \end{smallmatrix}\right], \left[\begin{smallmatrix} T' & 0 \\ 
M & U'  \end{smallmatrix}\right]  \big) := \!\!\left[\begin{smallmatrix}
\Ho_{\mathcal{T}}^{n}(T,T') & 0\\
M^{n}(U',T) & \Ho_{\mathcal{U}}^{n}(U,U')
\end{smallmatrix}\right]=\Ho_{\mathcal{T}}(T,T')^{n} \oplus M^{n}(U',T) $ $\oplus \Ho_{\mathcal{U}}^{n}(U,U')$  for all $n \in \Z$.\\
Now, for $\Delta:=\Ho_{\Lambda} \big( \left[\begin{smallmatrix} T & 0 \\ 
M & U  \end{smallmatrix}\right], \left[\begin{smallmatrix} T' & 0 \\ 
M & U'  \end{smallmatrix}\right]  \big)$ we define a differential
$$d_{\Delta}:\Ho_{\Lambda} \big( \left[\begin{smallmatrix} T & 0 \\ 
M & U  \end{smallmatrix}\right], \left[\begin{smallmatrix} T' & 0 \\ 
M & U'  \end{smallmatrix}\right]  \big)  \longrightarrow \Ho_{\Lambda} \big( \left[\begin{smallmatrix} T & 0 \\ 
M & U  \end{smallmatrix}\right], \left[\begin{smallmatrix} T' & 0 \\ 
M & U'  \end{smallmatrix}\right]  \big)  $$
$$ \hspace{2cm} \left[\begin{smallmatrix}
t & 0 \\
m & u
\end{smallmatrix}\right] \longmapsto \left[\begin{smallmatrix}
d_{\mathcal{T}(T,T')}(t) & 0 \\
d_{M(U',T)}(m) & d_{\mathcal{U}(U,U')}(u)
\end{smallmatrix}\right]. $$
 It is easy to see that $d_{\Delta}$ is a $K$-linear map of degree
$1$ and that $d_{\Delta}\circ d_{\Delta}=0$, then $d_{\Delta}$ is a differential for $\Delta$.
Then we have proved that 
$$\Ho_{\Lambda} \Big( {\left[ \begin{smallmatrix} T & 0 \\ 
M & U  \end{smallmatrix}\right]},{\left[\begin{smallmatrix} T' & 0 \\ 
M & U'  \end{smallmatrix}\right]}  \Big) := \left[\begin{smallmatrix}
{\Ho_{\mathcal{T}}(T,T')} & 0\\
M(U',T) &{ \Ho_{\mathcal{U}}(U,U')}
\end{smallmatrix}\right]$$
is a dg $K$-module.\\
Now let's see that the composition function
$$\Theta : \left[\begin{smallmatrix}
\mathcal{T}(T',T'') & 0\\
M(U'',T') & \mathcal{U}(U',U'')
\end{smallmatrix}\right] \times  \left[\begin{smallmatrix}
\mathcal{T}(T,T') & 0\\
M(U',T) & \mathcal{U}(U,U')
\end{smallmatrix}\right] \longrightarrow  \left[\begin{smallmatrix}
\mathcal{T}(T,T'') & 0\\
M(U'',T) & \mathcal{U}(U,U'')
\end{smallmatrix}\right] $$
$$\big( \left[\begin{smallmatrix}
t_{2} & 0 \\
m_{2} & u_{2}
\end{smallmatrix}\right],
\left[\begin{smallmatrix}
t_{1} & 0 \\
m_{1} & u_{1}
\end{smallmatrix}\right] \big)
\longmapsto 
\left[\begin{smallmatrix}
t_{2}\circ t_{1} & 0 \\
m_{2}\bullet t_{1}+ u_{2} \bullet m_{1} & u_{2}\circ u_{1}
\end{smallmatrix}\right],$$
induces  a morphism of dg-$K$-modules of zero degree.
In order to do that, let   $\left[\begin{smallmatrix}
T&0 \\
M&U
\end{smallmatrix}\right],$ $\left[\begin{smallmatrix}
T'&0 \\
M&U'
\end{smallmatrix}\right]$ and  $\left[\begin{smallmatrix}
T''&0 \\
M&U''
\end{smallmatrix}\right]
  \in \Lambda $ be; and for simplicity let us denote 
$$\Delta_{1}:=\Ho_{\Lambda} \big( \left[\begin{smallmatrix}
T & 0 \\
M & U
\end{smallmatrix}\right], \left[\begin{smallmatrix}
T' & 0 \\
M & U'
\end{smallmatrix}\right]\big),$$ 
$$\Delta_{2}:=\Ho_{\Lambda} \big( \left[\begin{smallmatrix}
T' & 0 \\
M & U'
\end{smallmatrix}\right], \left[\begin{smallmatrix}
T'' & 0 \\
M & U''
\end{smallmatrix}\right] \big),$$

$$\Delta_{3}:=\Ho_{\Lambda} \big( \left[\begin{smallmatrix}
T & 0 \\
M & U
\end{smallmatrix}\right], \left[\begin{smallmatrix}
T'' & 0 \\
M & U''
\end{smallmatrix}\right] \big).$$
We claim that the following diagram commutes 
$$\xymatrix{ { \Ho_{\Lambda} \big(    \left[\begin{smallmatrix}
T' & 0 \\
M & U'
\end{smallmatrix}\right],
\left[\begin{smallmatrix}
T'' & 0\\
M & U''
\end{smallmatrix}\right] \big) \otimes \Ho_{\Lambda} \big(    \left[\begin{smallmatrix}
T & 0 \\
M & U
\end{smallmatrix}\right],
\left[\begin{smallmatrix}
T' & 0\\
M & U'
\end{smallmatrix}\right] \big) }\ar[dd]_{D} \ar[r]^(.65){\Theta}&{\Ho_{\Lambda} \big(    \left[\begin{smallmatrix}
T & 0 \\
M & U
\end{smallmatrix}\right],
\left[\begin{smallmatrix}
T'' & 0\\
M & U''
\end{smallmatrix}\right] \big)} \ar[dd]^{d_{\Delta_{3}}}& \\
& & \\
{ \Ho_{\Lambda} \big(    \left[\begin{smallmatrix}
T' & 0 \\
M & U'
\end{smallmatrix}\right],
\left[\begin{smallmatrix}
T'' & 0\\
M & U''
\end{smallmatrix}\right] \big) \otimes \Ho_{\Lambda} \big(    \left[\begin{smallmatrix}
T & 0 \\
M & U
\end{smallmatrix}\right],
\left[\begin{smallmatrix}
T' & 0\\
M & U'
\end{smallmatrix}\right] \big) } \ar[r]_(.65){\Theta}&{\Ho_{\Lambda} \big(    \left[\begin{smallmatrix}
T & 0 \\
M & U
\end{smallmatrix}\right],
\left[\begin{smallmatrix}
T'' & 0\\
M & U''
\end{smallmatrix}\right] \big)}, &}$$
where $D$ is defined as in Equation \ref{tendgestructure}. \\
Indeed, for $\left[\begin{smallmatrix}
t_{2} & 0\\
m_{2} & u_{2}
\end{smallmatrix}\right] \otimes \left[\begin{smallmatrix}
t_{1} & 0\\
m_{1} & u_{1}
\end{smallmatrix}\right] \in \Ho_{\Lambda} \big( \left[\begin{smallmatrix}
T' & 0 \\
M & U'
\end{smallmatrix}\right], \left[\begin{smallmatrix}
T'' & 0 \\
M & U''
\end{smallmatrix}\right] \big) \otimes \Ho_{\Lambda} \big( \left[\begin{smallmatrix}
T & 0 \\
M & U
\end{smallmatrix}\right], \left[\begin{smallmatrix}
T' & 0 \\
M & U'
\end{smallmatrix}\right] \big)$
we get that

\begin{align*}
& D\Big(\left[\begin{smallmatrix}
t_{2} & 0\\
m_{2} & u_{2}
\end{smallmatrix}\right] \otimes \left[\begin{smallmatrix}
t_{1} & 0\\
m_{1} & u_{1}
\end{smallmatrix}\right]\Big)=\\
& = d_{\Delta_{2}}\Big(\left[\begin{smallmatrix}
t_{2} & 0\\
m_{2} & u_{2}
\end{smallmatrix}\right] \Big)\otimes \left[\begin{smallmatrix}
t_{1} & 0\\
m_{1} & u_{1}
\end{smallmatrix}\right] +(-1)^{n}\left[\begin{smallmatrix}
t_{2} & 0\\
m_{2} & u_{2}
\end{smallmatrix}\right]\otimes d_{\Delta_{1}}\Big(\left[\begin{smallmatrix}
t_{1} & 0\\
m_{1} & u_{1}
\end{smallmatrix}\right]\Big)\\
& =\left[\begin{smallmatrix}
d_{\mathcal{T}(T',T'')}(t_{2}) & 0\\
d_{M(U'',T')}(m_{2}) & d_{\mathcal{U}(U',U'')}(u_{2})
\end{smallmatrix}\right] \otimes \left[\begin{smallmatrix}
t_{1} & 0\\
m_{1} & u_{1}
\end{smallmatrix}\right] +(-1)^{n}\left[\begin{smallmatrix}
t_{2} & 0\\
m_{2} & u_{2}
\end{smallmatrix}\right]\otimes \left[\begin{smallmatrix}
d_{\mathcal{T}(T,T')}(t_{1}) & 0\\
d_{M(U',T)}(m_{1}) & d_{\mathcal{U}(U,U')}(u_{1})
\end{smallmatrix}\right]
\end{align*}
where $n:=|\left[\begin{smallmatrix}
t_{2} & 0\\
m_{2} & u_{2}
\end{smallmatrix}\right]|=|t_{2}|=|m_{2}|=|u_{2}|$.\\
 Hence, we obtain that:
 
\begin{align*}
& \Theta\Big(D\Big(\left[\begin{smallmatrix}
t_{2} & 0\\
m_{2} & u_{2}
\end{smallmatrix}\right] \otimes \left[\begin{smallmatrix}
t_{1} & 0\\
m_{1} & u_{1}
\end{smallmatrix}\right]\Big)\Big)\\
& =\left[\begin{smallmatrix}
d_{\mathcal{T}(T',T'')}(t_{2}) & 0\\
d_{M(U'',T')}(m_{2}) & d_{\mathcal{U}(U',U'')}(u_{2})
\end{smallmatrix}\right] \circ \left[\begin{smallmatrix}
t_{1} & 0\\
m_{1} & u_{1}
\end{smallmatrix}\right] +(-1)^{n}\left[\begin{smallmatrix}
t_{2} & 0\\
m_{2} & u_{2}
\end{smallmatrix}\right]\circ \left[\begin{smallmatrix}
d_{\mathcal{T}(T,T')}(t_{1}) & 0\\
d_{M(U',T)}(m_{1}) & d_{\mathcal{U}(U,U')}(u_{1})
\end{smallmatrix}\right]\\
& =\left[\begin{smallmatrix}
d_{\mathcal{T}(T',T'')}(t_{2})\circ t_{1} & 0\\
X, &\quad  d_{\mathcal{U}(U',U'')}(u_{2})\circ u_{1}
\end{smallmatrix}\right]+\left[\begin{smallmatrix}
(-1)^{n}t_{2}\circ d_{\mathcal{T}(T,T')}(t_{1}) & 0\\
Y, &\quad  (-1)^{n}u_{2}\circ d_{\mathcal{U}(U,U')}(u_{1})
\end{smallmatrix}\right]\\
& =  \left[\begin{smallmatrix} 
d_{\mathcal{T}(T',T'')}(t_{2})\circ t_{1}+(-1)^{n}t_{2}\circ d_{\mathcal{T}(T,T')}(t_{1})  & 0\\
X+Y, &  d_{\mathcal{U}(U',U'')}(u_{2})\circ u_{1}+  (-1)^{n}u_{2}\circ d_{\mathcal{U}(U,U')}(u_{1})
\end{smallmatrix}\right]\\
& =  \left[\begin{smallmatrix} 
d_{\mathcal{T}(T,T'')}(t_{2}\circ t_{1})   & 0\\
X+Y, &  d_{\mathcal{U}(U,U'')}(u_{2}\circ u_{1}) 
\end{smallmatrix}\right],
\end{align*}
since $n=|u_{2}|=|t_{2}|$, where $X:=d_{M(U'',T')}(m_{2})\bullet t_{1}+ d_{\mathcal{U}(U',U'')}(u_{2})\bullet m_{1}$ and 
$Y:=(-1)^{n}m_{2}\bullet d_{\mathcal{T}(T,T')}(t_{1}) + (-1)^{n}u_{2}\bullet d_{M(U',T)}(m_{1})$.\\
Now, we obtain that
\begin{align*}
X+Y& =\Big(d_{M(U'',T')}(m_{2})\bullet t_{1}+ d_{\mathcal{U}(U',U'')}(u_{2})\bullet m_{1}\Big)+\\
 & +\Big((-1)^{n}m_{2}\bullet d_{\mathcal{T}(T,T')}(t_{1}) + (-1)^{n}u_{2}\bullet d_{M(U',T)}(m_{1})\Big)\\
& =\Big(d_{M(U'',T')}(m_{2})\bullet t_{1}+(-1)^{n}m_{2}\bullet d_{\mathcal{T}(T,T')}(t_{1})\Big)\\
& +\Big(d_{\mathcal{U}(U',U'')}(u_{2})\bullet m_{1}+(-1)^{n}u_{2}\bullet d_{M(U',T)}(m_{1})\Big)=\\
& =d_{M(U'',T)}(m_{2}\bullet t_{1})+d_{M(U'',T)}(u_{2}\bullet m_{1})
\end{align*}
where the last equality is by Lemma \ref{similarreglaLeibniz}. We conclude that

\begin{align*}
& \Theta\Big(D\Big(\left[\begin{smallmatrix}
t_{2} & 0\\
m_{2} & u_{2}
\end{smallmatrix}\right] \otimes \left[\begin{smallmatrix}
t_{1} & 0\\
m_{1} & u_{1}
\end{smallmatrix}\right]\Big)\Big)=\left[\begin{smallmatrix} 
d_{\mathcal{T}(T,T'')}(t_{2}\circ t_{1})   & 0\\
d_{M(U'',T)}(m_{2}\bullet t_{1})+d_{M(U'',T)}(u_{2}\bullet m_{1}), & \quad d_{\mathcal{U}(U,U'')}(u_{2}\circ u_{1}) 
\end{smallmatrix}\right].
\end{align*}

On the other hand, we have the following equalities:

\begin{align*}
& d_{\Delta_{3}}\Big(\Theta\Big(\left[\begin{smallmatrix}
t_{2} & 0\\
m_{2} & u_{2}
\end{smallmatrix}\right] \otimes \left[\begin{smallmatrix}
t_{1} & 0\\
m_{1} & u_{1}
\end{smallmatrix}\right]\Big)\Big) =
d_{\Delta_{3}}\Big(\left[\begin{smallmatrix} 
t_{2}\circ t_{1} & 0\\
m_{2}\bullet t_{1}+u_{2}\bullet m_{1} & u_{2}\circ u_{1}
\end{smallmatrix}\right]\Big)=\\
& =\left[\begin{smallmatrix} 
d_{\mathcal{T}(T,T'')}(t_{2}\circ t_{1})  &\quad  0\\
d_{M(U'',T)}(m_{2}\bullet t_{1}+u_{2}\bullet m_{1})   &\quad d_{\mathcal{U}(U,U'')}(u_{2}\circ u_{1})
\end{smallmatrix}\right]\\
& =\left[\begin{smallmatrix} 
d_{\mathcal{T}(T,T'')}(t_{2}\circ t_{1})   & 0\\
d_{M(U'',T)}(m_{2}\bullet t_{1})+d_{M(U'',T)}(u_{2}\bullet m_{1}), & \quad d_{\mathcal{U}(U,U'')}(u_{2}\circ u_{1}) 
\end{smallmatrix}\right].
\end{align*}
This proves  that $d_{\Delta_{3}}\circ \Theta=\Theta\circ D,$ and hence $\Lambda=\left[\begin{smallmatrix} 
\mathcal{T} & 0\\
M & \mathcal{U}
\end{smallmatrix}\right]$ is a dg $K$-category.
\end{proof}

Now, consider an object $(A,f,B) \in \big( \mathrm{DgMod}(\T), \mathbb{G}\mathrm{DgMod}(\U)\big)$ in the dg-comma category. That is, $f\in \mathrm{Hom}_{\mathrm{DgMod}(\T)}^{0}(A,\mathbb{G}(B))$ and $d_{\mathrm{Hom}_{\mathrm{DgMod}(\mathcal{T})}(A,\mathbb{G}(B))}(f)=0$.  We can construct a functor  $A \coprod_{f} B: \Lambda \longrightarrow \mathrm{DgMod}(K)$ as follows .

\begin{enumerate}
\item[(a)] For  $\left[\begin{smallmatrix}
T & 0\\
M & U
\end{smallmatrix}\right] \in \Lambda $  we define:  $\big(  A \coprod_{f} B\big)\big(  \left[\begin{smallmatrix}
T & 0\\
M & U
\end{smallmatrix}\right]\big):= A(T) \coprod B(U) \in \mathrm{DgMod}(K)$ with differential
$d_{\big(  A \coprod_{f} B\big)\big(  \left[\begin{smallmatrix}
T & 0\\
M & U
\end{smallmatrix}\right]\big)}:=d_{A(T)}\coprod d_{B(U)}$.

\item[(b)] For an arbitrary homogenous element  $ \left[\begin{smallmatrix}
t & 0\\
m & u
\end{smallmatrix}\right] \in$ $\Ho_{\Lambda} \big(  \left[\begin{smallmatrix}
T & 0\\
M & U
\end{smallmatrix}\right],  \left[\begin{smallmatrix}
T' & 0\\
M & U'
\end{smallmatrix}\right] \big)=$ $ \left[\begin{smallmatrix}
\Ho_{\T}(T,T') & 0\\
M(U',T) & \Ho_{\U}(U,U')
\end{smallmatrix}\right]$  we define the morphism: 
$$\big( A \coprod_{f} B \big) \big(  \left[\begin{smallmatrix}
t & 0\\
m & u
\end{smallmatrix}\right] \big):=  \left[\begin{smallmatrix}
A(t) & 0\\
m & B(u)
\end{smallmatrix}\right]: A(T) \coprod B(U) \longrightarrow A(T') \coprod B(U')$$
 given by $ \left[\begin{smallmatrix}
A(t) & 0\\
m & B(u)
\end{smallmatrix}\right] \left[\begin{smallmatrix}
x\\
y
\end{smallmatrix}\right]=  \left[\begin{smallmatrix}
A(t)(x)\\
m \cdot x + B(u)(y)
\end{smallmatrix}\right]$ for each homogenous element $(x,y) \in A(T) \coprod B(U)$ (in this case, $|(x,y)|=|x|=|y|$), where $m \cdot x = (-1)^{|m||x|} [f_{T}(x)]_{U'}(m) \in B(U')$ (see Proposition \ref{productos dg matrices}).
\end{enumerate}
\begin{note}
In terms of  Proposition \ref{productos dg matrices},  we have that
$$ \left[\begin{smallmatrix}
A(t) & 0\\
m & B(u)
\end{smallmatrix}\right] \left[\begin{smallmatrix}
x\\
y
\end{smallmatrix}\right]=  \left[\begin{smallmatrix}
t \ast x\\
m \cdot x + u \diamond y
\end{smallmatrix}\right]\,\text{for each homogenous element }\,\,(x,y)\in A(T) \coprod B(U).$$
\end{note}
In order to prove that $A \coprod_{f} B: \Lambda \longrightarrow \mathrm{DgMod}(K)$ is a dg-functor, we need the following Lemma.

\begin{lemma}\label{otroparecidoLeinniz}
For homogenous elements $m\in M(U,T)$, $x\in A(T)$  we obtain that:
$$d_{B(U)}\big(m\cdot x\big)=d_{M(U,T)}(m)\cdot x+(-1)^{|m|}m\cdot d_{A(T)}(x).$$
\end{lemma}
\begin{proof}
We have that $f\in \mathrm{Hom}_{\mathrm{DgMod}(\T)}^{0}(A,\mathbb{G}(B))$ and
$D(f)=0$ where 
$$D:=d_{\mathrm{Hom}_{\mathrm{DgMod}(\mathcal{T})}(A,\mathbb{G}(B))}:\mathrm{Hom}_{\mathrm{DgMod}(\T)}(A,\mathbb{G}(B))\longrightarrow \mathrm{Hom}_{\mathrm{DgMod}(\T)}(A,\mathbb{G}(B)).$$
Recall that for $T\in \mathcal{T}$ we get that 
$D(f)_{T}:=D'(f_{T})$ (see Equation \ref{dgNatmod})  where 
$$D':\mathrm{Hom}_{\mathrm{DgMod}(K)}\big(A(T),\mathbb{G}(B)(T)\big)\longrightarrow \mathrm{Hom}_{\mathrm{DgMod}(K)}\big(A(T),\mathbb{G}(B)(T)\big).$$
Since $D(f)=0$ we have that $D(f)_{T}:=D'(f_{T})=0$ for all $T\in \mathcal{T}$.
By Equation \ref{Homdgestructure}, it follows that $d_{\mathbb{G}(B)(T)}\circ f_{T}=f_{T}\circ d_{A(T)}$.
Then, on one hand, for $x\in A(T)$ homogenous, we have that $f_{T}(x)\in \mathrm{Hom}_{\mathrm{DgMod}(\U)}^{|x|}(M_{T},B)$. Hence $d_{\mathbb{G}(B)(T)}(f_{T}(x))\in \mathrm{Hom}_{\mathrm{DgMod}(\U)}^{|x|+1}(M_{T},B)$ where by Equation \ref{dgNatmod} obtain that:
\begin{align*}
[d_{\mathbb{G}(B)(T)}(f_{T}(x))]_{U}:=D''([f_{T}(x)]_{U})
\end{align*}
with $D'':\mathrm{Hom}_{\mathrm{DgMod}(K)}(M(U,T),B(U))\longrightarrow \mathrm{Hom}_{\mathrm{DgMod}(K)}(M(U,T),B(U))$.\\
By Equation \ref{Homdgestructure}, it follows that:
\begin{align*}
D''([f_{T}(x)]_{U}) & =d_{B(U)}\circ [f_{T}(x)]_{U}-(-1)^{[f_{T}(x)]_{U}}[f_{T}(x)]_{U}\circ d_{M(U,T)}\\
 & =d_{B(U)}\circ [f_{T}(x)]_{U}-(-1)^{|x|}[f_{T}(x)]_{U}\circ d_{M(U,T)}.
\end{align*}

On the other hand, $(f_{T}\circ d_{A(T)})(x)=f_{T}(d_{A(T)}(x))\in \mathbb{G}(B)(T)$, where for $U\in \mathcal{U}$ we have $[f_{T}(d_{A(T)}(x))]_{U}:M(U,T)\longrightarrow B(U)$. Hence we have proved that
 
 \begin{align*}
 [f_{T}(d_{A(T)}(x))]_{U}=d_{B(U)}\circ [f_{T}(x)]_{U}-(-1)^{|x|}[f_{T}(x)]_{U}\circ d_{M(U,T)}.
 \end{align*}
 Then for $m\in M(U,T)$ homogenous we have that
 \begin{align*}
 [f_{T}(d_{A(T)}(x))]_{U}(m)=\big(d_{B(U)}\circ [f_{T}(x)]_{U}\big)(m)-(-1)^{|x|}\big([f_{T}(x)]_{U}\circ d_{M(U,T)}\big)(m).
 \end{align*}
 Hence,
\begin{align*}
& m\cdot d_{A(T)}(x) =(-1)^{|m|\cdot |d_{A(T)}(x)|}[f_{T}(d_{A(T)}(x))]_{U}(m)\\
& =(-1)^{|m|\cdot (|x|+1)}[f_{T}(d_{A(T)}(x))]_{U}(m)\\
& = (-1)^{|m|}(-1)^{|m||x|}\big(d_{B(U)}\circ [f_{T}(x)]_{U}\big)(m) +\\
& -(-1)^{|x|} (-1)^{|m|}(-1)^{|m||x|}\big([f_{T}(x)]_{U}\circ d_{M(U,T)}\big)(m)\\
& = (-1)^{|m|}d_{B(U)}\big(m\cdot x\big)-(-1)^{|m|} (-1)^{(|m|+1)|x|}\big([f_{T}(x)]_{U}\circ d_{M(U,T)}\big)(m)\\
& = (-1)^{|m|}d_{B(U)}\big(m\cdot x\big)-(-1)^{|m|} d_{M(U,T)}(m)\cdot x.
\end{align*}
Hence we obtain that: $d_{B(U)}\big(m\cdot x\big)=d_{M(U,T)}(m)\cdot x+(-1)^{|m|}m\cdot d_{A(T)}(x).$
\end{proof}

The following lemma tell us that $A \coprod_{f} B$ is a  dg-functor.

\begin{lemma}\label{A copro B df-funtor}
Let $(A, f,B) \in  \big( \mathrm{DgMod}(\T), \mathbb{G}(\mathrm{DgMod}(\U)\big)$ be. Then  $A \coprod_{f} B: \Lambda \longrightarrow \mathrm{DgMod}(K)$ is a dg-functor.
\end{lemma}
\begin{proof}
Using Proposition \ref{productos dg matrices}(c) it  follows that $A \coprod_{f} B$ is a functor. Now let us check that $\big( A  \coprod_{f} B \big)$ is a graded functor, that is, $\big( A \coprod_{f} B \big)\Big(\Ho_{\Lambda}^{n} \big( \left[\begin{smallmatrix}
T & 0\\
M & U
\end{smallmatrix}\right],$\\
$ \left[\begin{smallmatrix}
T' & 0\\
M & U'
\end{smallmatrix}\right] \big)\Big)\subseteq \Ho_{\Dg(K)}^{n} \big(  \big(A \coprod_{f} B \big) \left[\begin{smallmatrix}
T & 0\\
M & U
\end{smallmatrix}\right], \big(A \coprod_{f} B \big) \left[\begin{smallmatrix}
T' & 0\\
M & U'
\end{smallmatrix}\right]\big).$\\
Indeed, let $\left[\begin{smallmatrix}
 t & 0 \\
 m & u
 \end{smallmatrix}\right] \in \Ho_{\Lambda}^{n} \big( \left[\begin{smallmatrix}
T & 0\\
M & U
\end{smallmatrix}\right], \left[\begin{smallmatrix}
T' & 0\\
M & U'
\end{smallmatrix}\right] \big)$ be a homogenous element, by definition we have that $\big(A \coprod_{f} B \big) \big( \left[\begin{smallmatrix}
t & 0 \\
 m & u
\end{smallmatrix}\right]\big)=\left[\begin{smallmatrix}
A(t) & 0 \\
 m & B(u)
\end{smallmatrix}\right] $. Since $t \in \Ho_{\T}^{n}(T,T')$, $u \in \Ho_{\U}^{n}(U,U')$, $m \in M^{n}(U',T)$ and $A,B$ are dg-functors it follows that $A(t) \in  \Ho_{\Dg(K)}^{n}(A(T),A(T'))$, $B(u) \in  \Ho_{\Dg(K)}^{n}(B(U),B(U'))$.\\
Now, let us see that 
$\left[\begin{smallmatrix}
A(t) & 0 \\
 m & B(u)
\end{smallmatrix}\right] \in   \Ho_{\Dg(K)}^{n} \big( (A \coprod_{f} B) \left[\begin{smallmatrix}
T & 0\\
M & U
\end{smallmatrix}\right], \big(A \coprod_{f} B \big) \left[\begin{smallmatrix}
T' & 0\\
M & U'
\end{smallmatrix}\right]\big).$\\ 
Indeed, let $(x,y)\in \big(  A \coprod_{f} B\big)\big(  \left[\begin{smallmatrix}
T & 0\\
M & U
\end{smallmatrix}\right]\big):= A(T) \coprod B(U)$ a homogenous element with $k=|(x,y)|=|x|=|y|$.\\
Then  $ \left[\begin{smallmatrix}
A(t) & 0\\
m & B(u)
\end{smallmatrix}\right] \left[\begin{smallmatrix}
x\\
y
\end{smallmatrix}\right]=  \left[\begin{smallmatrix}
A(t)(x)\\
m \cdot x + B(u)(y)
\end{smallmatrix}\right]$  where $m \cdot x = (-1)^{|m||x|} [f_{T}(x)]_{U'}(m) \in B(U')$. We note that
$|A(t)(x)|=|t|+|x|=n+k$, $|B(u)(y)|=|u|+|y|=n+k$ and  by Proposition \ref{productos dg matrices}, we have that 
$|m \cdot x|=|m|+|x|=n+k$. Then $\left[\begin{smallmatrix}
A(t)(x)\\
m \cdot x + B(u)(y)
\end{smallmatrix}\right]$  has degree $n+k$.\\
This proves that $\left[\begin{smallmatrix}
A(t) & 0 \\
 m & B(u)
\end{smallmatrix}\right] \in   \Ho_{\Dg(K)}^{n} \big( (A \coprod_{f} B) \left[\begin{smallmatrix}
T & 0\\
M & U
\end{smallmatrix}\right], \big(A \coprod_{f} B \big) \left[\begin{smallmatrix}
T' & 0\\
M & U'
\end{smallmatrix}\right]\big).$ Therefore,  $\big( A  \coprod_{f} B \big)$ is a graded functor.\\
Now, let us see that $ A \coprod_{f} B$ commutes with the differentials. For simplicity we set $\Delta_{1}:=\Ho_{\Lambda} \big( \left[\begin{smallmatrix}
T & 0\\
M & U
\end{smallmatrix}\right], \left[\begin{smallmatrix}
T' & 0\\
M & U'
\end{smallmatrix}\right] \big)$. Let us check that the following diagram commutes
$$\xymatrix{{ \Ho_{\Lambda} \big( \left[\begin{smallmatrix}
T & 0\\
M & U
\end{smallmatrix}\right], \left[\begin{smallmatrix}
T' & 0\\
M & U'
\end{smallmatrix}\right] \big)} \ar[r]^(.37){(A \coprod_{f} B)} \ar[d]_{d_{\Delta_{1}}}  & \Ho_{\Dg(K)} \big( A(T)\coprod B(U), A(T')\coprod B(U') \big) \ar[d]^{d_{\Delta_{2}}} \\
{ \Ho_{\Lambda} \big( \left[\begin{smallmatrix}
T & 0\\
M & U
\end{smallmatrix}\right], \left[\begin{smallmatrix}
T' & 0\\
M & U'
\end{smallmatrix}\right] \big)}\ar[r]_(.37){(A\coprod_{f} B)} & \Ho_{\Dg(K)} \big( A(T)\coprod B(U), A(T')\coprod B(U') \big) }$$
for all $\left[\begin{smallmatrix}
T & 0\\
M & U
\end{smallmatrix}\right], \left[\begin{smallmatrix}
T' & 0\\
M & U'
\end{smallmatrix}\right] \in \Lambda $.\\

Indeed, let $\left[\begin{smallmatrix}
t & 0\\
m & u
\end{smallmatrix}\right] \in  \Ho_{\Lambda} \big( \left[\begin{smallmatrix}
T & 0\\
M & U
\end{smallmatrix}\right], \left[\begin{smallmatrix}
T' & 0\\
M & U'
\end{smallmatrix}\right] \big)$ be a homogenous element of degree $n$ (that is, $n=|t|=|m|=|u|$). On one hand, we have that

\begin{align*}
& \Big(d_{\Delta_{2}}\circ (A \coprod_{f} B )\Big) \big( \left[\begin{smallmatrix}
t & 0\\
m & u
\end{smallmatrix}\right] \big)  = d_{\Delta_{2}} \Big( \left[\begin{smallmatrix}
A(t) & 0\\
m & B(u)
\end{smallmatrix}\right] \Big)=\\
&=(d_{A(T')}\coprod d_{B(U')})\circ  \left[\begin{smallmatrix}
A(t) & 0\\
m & B(u)
\end{smallmatrix}\right]- (-1)^{n}\left[\begin{smallmatrix}
A(t) & 0\\
m & B(u)
\end{smallmatrix}\right] \circ (d_{A(T)}\coprod d_{B(U)}).
\end{align*}
Then for $(x,y)\in A(T)\coprod B(U)$ homogenous element of degree $|(x,y)|=|x|=|y|=k$ we get that:

\begin{align*}
&=\Big((d_{A(T')}\coprod d_{B(U')})\circ  \left[\begin{smallmatrix}
A(t) & 0\\
m & B(u)
\end{smallmatrix}\right]- (-1)^{n}\left[\begin{smallmatrix}
A(t) & 0\\
m & B(u)
\end{smallmatrix}\right] \circ (d_{A(T)}\coprod d_{B(U)})\Big)(x,y)=\\
&=(d_{A(T')}\coprod d_{B(U')})\Big(  \left[\begin{smallmatrix}
A(t)(x) \\
m\cdot x +B(u)(y)
\end{smallmatrix}\right]\Big)- (-1)^{n}\left[\begin{smallmatrix}
A(t) & 0\\
m & B(u)
\end{smallmatrix}\right]  \left[\begin{smallmatrix}
d_{A(T)}(x) \\
d_{B(U)}(y)
\end{smallmatrix}\right]\\
&=\left[\begin{smallmatrix}
d_{A(T')}(A(t)(x)) \\
d_{B(U')}\big(m\cdot x +B(u)(y)\big)
\end{smallmatrix}\right]- (-1)^{n}\left[\begin{smallmatrix}
A(t)(d_{A(T)}(x)) \\
m\cdot d_{A(T)}(x) + B(u)(d_{B(U)}(y))
\end{smallmatrix}\right]\\
& =\left[\begin{smallmatrix}
d_{A(T')}(A(t)(x)) -(-1)^{n}A(t)(d_{A(T)}(x))  \\
d_{B(U')}\big(m\cdot x +B(u)(y)\big)-(-1)^{n}m\cdot d_{A(T)}(x)  -(-1)^{n}B(u)(d_{B(U)}(y))
\end{smallmatrix}\right]\\
& =\left[\begin{smallmatrix}
d_{A(T')}(A(t)(x)) -(-1)^{n}A(t)(d_{A(T)}(x))  \\
d_{B(U')}\big(m\cdot x\big) +d_{B(U')}\big( B(u)(y)\big)-(-1)^{n}m\cdot d_{A(T)}(x) -(-1)^{n}B(u)(d_{B(U)}(y))
\end{smallmatrix}\right]\\
& =\left[\begin{smallmatrix}
d_{A(T')}(A(t)(x)) -(-1)^{n}A(t)(d_{A(T)}(x))  \\
\Big(d_{B(U')}\big(m\cdot x\big) -(-1)^{n}m\cdot d_{A(T)}(x)\Big)+\Big(d_{B(U')}\big( B(u)(y)\big) -(-1)^{n}B(u)(d_{B(U)}(y))\Big)
\end{smallmatrix}\right]
\end{align*}

Now, since $A:\mathcal{T}\longrightarrow \mathrm{DgMod}(K)$ is a dg-functor,  for $x\in A(t)$ homogenous of degree $k$ and $t$ homogenous of degree $n$ we have that
\begin{align*}
A(d_{\mathcal{T}(T,T')}(t))(x)= d_{A(T')}(A(t)(x))-(-1)^{n}A(t)(d_{A(T)}(x)).
\end{align*}

Similarly, since $B:\mathcal{U}\longrightarrow \mathrm{DgMod}(K)$ is a dg-functor, for $y\in B(U)$ homogenous of degree $k$ and $u$ homogenous of degree $n$ we have that
\begin{align*}
B(d_{\mathcal{U}(U,U')}(u))(y)= d_{B(U')}(B(u)(y))-(-1)^{n}B(u)(d_{B(U)}(y)).
\end{align*}
Since $m\in M(U',T)$ and $x\in A(t)$ with $|m|=n$ and $|x|=k$, by Lemma \ref{otroparecidoLeinniz}, we obtain that 
$d_{B(U')}\big(m\cdot x\big)=d_{M(U',T)}(m)\cdot x+(-1)^{|m|}m\cdot d_{A(T)}(x).$ Therefore, we conclude that $
\Big(d_{\Delta_{2}}\circ (A \coprod_{f} B )\Big) \big( \left[\begin{smallmatrix}
t & 0\\
m & u
\end{smallmatrix}\right] \big)(x,y)=\left[\begin{smallmatrix}
A(d_{\mathcal{T}(T,T')}(t))(x)\\
d_{M(U',T)}(m)\cdot x+B(d_{\mathcal{U}(U,U')}(u))(y)
\end{smallmatrix}\right].$\\
On the other hand,
\begin{align*}
(A\coprod_{f}B)\Big(d_{\Delta_{1}}\big( \left[\begin{smallmatrix}
t & 0\\
m & u
\end{smallmatrix}\right]\big)\Big) & =(A\coprod_{f}B)\Big(\left[\begin{smallmatrix}
d_{\mathcal{T}(T,T')}(t) & 0\\
d_{M(U',T)}(m) & d_{\mathcal{U}(U,U')}(u)
\end{smallmatrix}\right]\Big)\\
& =\left[\begin{smallmatrix}
A(d_{\mathcal{T}(T,T')}(t)) & 0\\
d_{M(U',T)}(m) & B(d_{\mathcal{U}(U,U')}(u))
\end{smallmatrix}\right]
\end{align*}

Then for $(x,y)\in A(T)\coprod B(U)$ homogenous element of degree $|(x,y)|=|x|=|y|=k$ we have that:

\begin{align*}
\left[\begin{smallmatrix}
A(d_{\mathcal{T}(T,T')}(t)) & 0\\
d_{M(U',T)}(m) & B(d_{\mathcal{U}(U,U')}(u))
\end{smallmatrix}\right](x,y)=\left[\begin{smallmatrix}
A(d_{\mathcal{T}(T,T')}(t))(x)\\
d_{M(U',T)}(m)\cdot x+B(d_{\mathcal{U}(U,U')}(u))(y)
\end{smallmatrix}\right].
\end{align*}
Hence, $\Big(d_{\Delta_{2}}\circ (A \coprod_{f} B )\Big) \big( \left[\begin{smallmatrix}
t & 0\\
m & u
\end{smallmatrix}\right] \big)(x,y)=(A\coprod_{f}B)\Big(d_{\Delta_{1}}\big( \left[\begin{smallmatrix}
t & 0\\
m & u
\end{smallmatrix}\right]\big)\Big)(x,y).$ Thus we get that 
$d_{\Delta_{2}}\circ (A \coprod_{f} B )=(A\coprod_{f}B)\circ d_{\Delta_{1}}.$ This proves that $A \coprod_{f} B :\Lambda\longrightarrow \mathrm{DgMod}(K)$ is a dg-functor.

\end{proof}

In this way we can construct a functor: $\mathfrak{F} : \Big( \Dg(\T), \mathbb{G}\big(\Dg(\U) \big)\Big) \longrightarrow \Dg(\Lambda)$
which is defined as follows.
\begin{itemize}
\item[(a)] For $(A,f,B) \in \Big(\Dg(\T), \mathbb{G}\big(\Dg(\U)\big)\Big)$  we define $\mathfrak{F}((A,f,B)):= A \coprod_{f} B.$
\item[(b)] If $(\alpha, \beta):(A,f,B) \longrightarrow (A',f',B')$ is a homogenous morphism of degree $n$ in $\Big(\Dg(\T), \mathbb{G}\Dg(\U)\Big)$ then $\mathfrak{F}(\alpha, \beta):= \alpha \coprod \beta$ is a  dg-natural transformation whose components are: 
$$\left \{\!\! (\alpha \coprod \beta)_{\left[\begin{smallmatrix}
T & 0\\
M & U\end{smallmatrix}\right]} := \alpha_{T} \coprod \beta_{U}: (A \coprod_{f} B) \big( \left[\begin{smallmatrix}
T & 0\\
M & U\end{smallmatrix}\right] \big) \longrightarrow (A' \coprod_{f'} B') \big( \left[\begin{smallmatrix}
T  & 0\\
M & U\end{smallmatrix}\right] \big)\!\! \right\}_{ \left[\begin{smallmatrix}
T & 0\\
M & U\end{smallmatrix}\right] \in \Lambda}$$
\end{itemize}

In the following Lemma we will see that $\mathfrak{F}(\alpha, \beta):= \alpha \coprod \beta$ is indeed a  dg-natural transformation. First, we recall that by the construction of the dg-comma category (see Proposition  \ref{dgcommaprop}), we have that if $(\alpha, \beta):(A,f,B) \longrightarrow (A',f',B')$ is a homogenous morphism of degree $n$ in $\Big(\Dg(\T), \mathbb{G}\Dg(\U)\Big)$, then 
$|\alpha|=|\beta|=n$.

  \begin{lemma}\label{a copro b  dg transformacion}
Let $(\alpha, \beta): (A,f,B) \longrightarrow (A',f',B')$ be a morphism in the comma category $\Big(\Dg(\T), \mathbb{G}\big(\Dg(\U)\big)\Big)$ of degree $n$ (that is, $n=|(\alpha,\beta)|= |\alpha|=|\beta|$). Then $\alpha \coprod \beta$ is a natural dg-transformation of degree $n=|\alpha|=|\beta|$.
\end{lemma}
\begin{proof}
Let $ \left[\begin{smallmatrix}
t_{1} & 0\\
m_{1} & u_{1}\end{smallmatrix}\right]:  \left[\begin{smallmatrix}
T & 0\\
M & U\end{smallmatrix}\right] \longrightarrow  \left[\begin{smallmatrix}
T' & 0\\
M & U'\end{smallmatrix}\right]$ be a  homogeneous morphism of degree  $r=|t_{1}|=|u_{1}|=|m_{1}|$ in $\Lambda= \left[ \begin{smallmatrix}
\T & 0 \\
M & \U
 \end{smallmatrix}\right]$.  We have to check   that the following  diagram commutes in $\Dg(K)$ up to the sign $(-1)^{nr}$ (see Definition \ref{dgnaturaltrans}):
$$(\ast):\xymatrix{ A(T) \coprod B(U) \ar[rr]^{\alpha_{T} \coprod \beta_{U}} \ar[dd]_{\left[\begin{smallmatrix}
A(t_{1}) & 0\\
m_{1} & B(u_{1})\end{smallmatrix}\right]}& &A'(T) \coprod B'(U)\ar[dd]^{ \left[\begin{smallmatrix}
A'(t_{1}) & 0\\
m_{1} & B'(u_{1})\end{smallmatrix}\right]} \\
& & \\
A(T') \coprod B(U') \ar[rr]_{\alpha_{T'} \coprod \beta_{U'}}& &A'(T') \coprod B'(U'). }$$

Indeed, for $(x,y) \in A(T) \coprod B(U) $ homogenous element of degree $|(x,y)|=|x|=|y|=k$, we have that:
\begin{align*}
\big( \left[ \begin{smallmatrix}
A'(t_{1}) & 0\\
m_{1} & B'(u_{1})
 \end{smallmatrix}\right] \circ (\alpha_{T} \coprod \beta_{U})\big) \left[ \begin{smallmatrix}
x\\
y
 \end{smallmatrix}\right]& = \left[ \begin{smallmatrix}
A'(t_{1}) & 0\\
m_{1} & B'(u_{1})
 \end{smallmatrix}\right] \left[ \begin{smallmatrix}
\alpha_{T}(x)\\
\beta_{U}(y)
 \end{smallmatrix}\right]\\
& = \left[ \begin{smallmatrix}
A'(t_{1})\big(\alpha_{T}(x)\big)\\
m_{1}\cdot \alpha_{T}(x)+ B'(u_{1})(\beta_{U}(y))
 \end{smallmatrix}\right].
\end{align*}
On the other hand we get that:
\begin{align*}
\big( (\alpha_{T'} \coprod \beta_{U'})\circ \left[ \begin{smallmatrix}
A(t_{1}) & 0\\
m_{1} & B(u_{1})
 \end{smallmatrix}\right] \big) \left[ \begin{smallmatrix}
x\\
y
 \end{smallmatrix}\right] & =(\alpha_{T'} \coprod \beta_{U'}) \left[ \begin{smallmatrix}
A(t_{1})(x)\\
m_{1} \cdot x + B(u_{1})(y)
 \end{smallmatrix}\right] \\
& = \left[ \begin{smallmatrix}
\alpha_{T'}(A(t_{1})(x))\\
\beta_{U'}\big( m_{1} \cdot x + B(u_{1})(y) \big)
 \end{smallmatrix}\right]\\
& =\left[ \begin{smallmatrix}
\alpha_{T'}(A(t_{1})(x))\\
\beta_{U'}( m_{1} \cdot x )+\beta_{U'} (B(u_{1})(y)) 
 \end{smallmatrix}\right].
\end{align*}

We will check that:

\begin{align*}
\big( \left[ \begin{smallmatrix}
A'(t_{1}) & 0\\
m_{1} & B'(u_{1})
 \end{smallmatrix}\right] \circ (\alpha_{T} \coprod \beta_{U})\big) \left[ \begin{smallmatrix}
x\\
y
 \end{smallmatrix}\right] =(-1)^{nr}\big( (\alpha_{T'} \coprod \beta_{U'})\circ \left[ \begin{smallmatrix}
A(t_{1}) & 0\\
m_{1} & B(u_{1})
 \end{smallmatrix}\right] \big) \left[ \begin{smallmatrix}
x\\
y
 \end{smallmatrix}\right]
\end{align*}
where $n=|(\alpha,\beta)|=|\alpha|=|\beta|$ and  $r=|t_{1}|=|u_{1}|=|m_{1}|$.\\
Firstly, since  $\alpha: A \longrightarrow A'$ is  a morphism in $\Dg(\T)$  of degree $n:=|\alpha|$, we have that $|\alpha_{T}|=n$ for all $T\in \mathcal{T}$. Then for $t_{1}: T \longrightarrow T'$ in $\T$ of degree $r=|t_{1}|$, we get that $A'(t_{1}) \alpha_{T}= (-1)^{|\alpha||t_{1}|} \alpha_{T'}A(t_{1})$.\\
 This implies that for an homogenous element $x \in A(T)$, we get that $A'(t_{1}) \Big(\alpha_{T}(x)\Big)= (-1)^{|\alpha||m_{1}|}\alpha_{T'}\Big(A(t_{1})(x)\Big)=  (-1)^{nr}\alpha_{T'}\Big(A(t_{1})(x)\Big)$.\\
 
Secondly, let's see that $m_{1}\cdot \alpha_{T}(x)= (-1)^{nr}\beta_{U'}( m_{1} \cdot x ).$\\
Indeed, for a homogenous elements $x \in A(T)$ and $m_{1}\in M(U',T)$ we have that

\begin{equation} \label{unaigual}
\begin{split}
 m_{1} \cdot \alpha_{T}(x)  = (-1)^{nr}(-1)^{|m_{1}||x|}[f'_{T}(\alpha_{T}(x))]_{U'}(m_{1}).
\end{split}
\end{equation}
 
Now, consider $\beta: B \longrightarrow B'$  in $\Dg(\U)$.  Then

\begin{equation} \label{dosigual}
\begin{split}
\beta_{U'}(m_{1} \cdot x) = \beta_{U'} \Big((-1)^{|m_{1}||x|} \big[f_{T}(x) \big]_{U'}(m_{1})\Big) =(-1)^{|m_{1}||x|}  \beta_{U'} \Big(\big[ f_{T}(x) \big]_{U'}(m_{1})\Big).
\end{split}
\end{equation}
 Since $(\alpha, \beta)$ is a morphism in the category $ \Big( \Dg(\T), \mathbb{G}\big(\Dg(\U)\big)\Big)$, for each $T \in \T$,  it follows that  $[f' \circ \alpha]_{T}= [ \mathbb{G}\beta \circ f]_{T}$. Then  $f'_{T}\circ \alpha_{T}= \Ho_{\Dg(\U)}(M_{T}, \beta)\circ  f_{T}$ (see  Note \ref{Obs de G}).  We note that $f_{T}(x)$ is a natural dg-transformation of degree $|x|$. Thus, for $m_{1} \in M(U',T)$  we have 
\begin{align*}
[f'_{T}\big(\alpha_{T}(x)\big)]_{U'}(m_{1})& = [\mathbb{G}(\beta)\big(f_{T}(x)\big)]_{U'}(m_{1})\\
&= [\Ho_{\Dg(\U)}(M_{T}, \beta) \big(f_{T}(x)\big)]_{U'}(m_{1})\\
& = [\beta \circ f_{T}(x)]_{U'}(m_{1})\\
& = \beta_{U'}\Big([f_{T}(x)]_{U'}(m_{1})\Big).
\end{align*}
Using the above equality together with (\ref{unaigual}) and (\ref{dosigual}) we get:
 \begin{align*}
 m_{1} \cdot \alpha_{T}(x) &= (-1)^{nr} \beta_{U'}(m_{1}\cdot x).
 \end{align*}
 
Finally, since $\beta: B \longrightarrow B'$ is a homogeneous morphism of degree $n=|\beta|$ in  $\Dg(\U)$, for  $u_{1}: U \longrightarrow U'$ homogenous of degree $r=|u_{1}|$ and a homogenous element $y \in B(U)$  we have that $B'(u_{1})\big(\beta_{U} (y)\big)= (-1)^{nr} \beta_{U'}\big( B(u_{1})(y)\big).$ 
All this proves that the diagram $(\ast)$ commute up to the sign $(-1)^{nr}$. Then  $\alpha \coprod \beta: A \coprod_{f} B \longrightarrow A' \coprod_{f'} B'$ is a natural dg-transformation of degree $|\alpha|=|\beta|.$
\end{proof}

\begin{proposition}\label{Fdgfunctor}
The assignment
$\mathfrak{F}: \Big( \Dg(\T), \mathbb{G}\big(\Dg(\U)\big) \Big) \longrightarrow \Dg(\Lambda) $
is a dg-functor.
\end{proposition}
\begin{proof}
Let us see that $\mathfrak{F}$ is graded, this is, that following  conditions is satisfied
$$\mathfrak{F} \big(  \Ho_{(\Dg(\T), \mathbb{G} \Dg(\U))}^{n}  ((A,f,B), (A',f',B'))\Big) \subseteq \Ho_{\Lambda}^{n} \big( A \amalg_{f}B, A' \amalg_{f} B' \big).$$
But, this is follows by  Lemma \ref{a copro b  dg transformacion}. Finally, let us check that $\mathfrak{F}$ commutes with the differential. Let $(A,f,B), (A',f',B') \in \mathcal{C}:=\Big(\Dg(\T), \mathbb{G}\big( \Dg(\U)\big)\Big) $ be and let us see that the following diagram commutes:
$$\xymatrix{ \Ho_{\mathcal{C}}((A,f,B),(A',f',B')) \ar[r]^(.5){\mathfrak{F}} \ar[d]_{d_{\Delta_{1}}}&  \Ho_{\Dg(\Lambda)}(A\amalg_{f}B, A' \amalg_{f'} B')\ar[d]^{d_{\Delta_{2}}} \\
\Ho_{\mathcal{C}}((A,f,B),(A',f',B')) \ar[r]_(.5){\mathfrak{F}}  & \Ho_{\Dg(\Lambda)}(A\amalg_{f}B, A' \amalg_{f'} B').}$$

Indeed, let $(\alpha, \beta) \in \Ho_{(\Dg(\T), \mathbb{G} \Dg(\U))}((A,f,B),(A',f',B')) $ be homogenous element of degree
$n=|(\alpha,\beta)|=|\alpha|=|\beta|$. Then we have that: 
$d_{\Delta_{2}}\mathfrak{F}((\alpha, \beta)) = d_{\Delta_{2}}(\alpha \amalg \beta).$ By Equation \ref{dgNatmod}, we have that for $\left[ \begin{smallmatrix}
T & 0\\
M & U
 \end{smallmatrix}\right] \in \Lambda$:
 
\begin{align*}
[d_{\Delta_{2}}(\alpha \amalg \beta)]_{\left[ \begin{smallmatrix}
T & 0\\
M & U
 \end{smallmatrix}\right] }:=D\Big([(\alpha \amalg \beta)]_{\left[ \begin{smallmatrix}
T & 0\\
M & U
 \end{smallmatrix}\right] }\Big)
\end{align*}

where $[(\alpha \amalg \beta)]_{\left[ \begin{smallmatrix}
T & 0\\
M & U
 \end{smallmatrix}\right]}=\alpha_{T}\coprod \beta_{U}:A(T)\coprod B(U)\longrightarrow A'(T)\coprod B'(U)$
 and $D$ is the differential of $\mathrm{Hom}_{\mathrm{DgMod}(K)}\Big(A(T)\coprod B(U),A'(T)\amalg B'(U)\Big)$. By Equation \ref{Homdgestructure}:
 
 \begin{align*}
& D\Big(\alpha_{T}\coprod \beta_{U}\Big) = \\
& =(d_{A'(T)}\coprod d_{B'(U)})\circ(\alpha_{T}\coprod \beta_{U})-(-1)^{r}(\alpha_{T}\coprod \beta_{U}) \circ (d_{A(T)}\coprod d_{B(U)})\\
 & =(d_{A'(T)}\circ \alpha_{T}\coprod d_{B'(U)}\circ \beta_{U}) -(-1)^{r}(\alpha_{T}\circ d_{A(T)}\coprod \beta_{U}\circ d_{B(U)})\\
& =\Big(d_{A'(T)}\circ \alpha_{T}-(-1)^{r}\alpha_{T}\circ d_{A(T)}\Big)\coprod \Big(d_{B'(U)}\circ \beta_{U}-(-1)^{r}\beta_{U}\circ d_{B(U)}\Big),
 \end{align*}
 where $r=|[(\alpha \amalg \beta)]_{\left[ \begin{smallmatrix}
T & 0\\
M & U
 \end{smallmatrix}\right]}|=|\alpha_{T}\coprod \beta_{U}|=n$.\\
On the other hand, by the definition of the differential in the comma category, we have that $d_{\Delta_{1}}(\alpha,\beta)=\Big(d_{\Ho_{\Dg(\T)} (A,A')}(\alpha),d_{\Ho_{\Dg(\U)} (B,B')}(\beta)\Big).$ Then 
\begin{align*}
 \mathfrak{F}\Big(d_{\Delta_{1}}(\alpha,\beta)\Big)=& \mathfrak{F}\Big(\Big(d_{\Ho_{\Dg(\T)} (A,A')}(\alpha),d_{\Ho_{\Dg(\U)} (B,B')}(\beta)\Big)\Big)\\
& = d_{\Ho_{\Dg(\T)} (A,A')}(\alpha)\coprod d_{\Ho_{\Dg(\U)} (B,B')}(\beta).
\end{align*}
Hence, for $\left[ \begin{smallmatrix}
T & 0\\
M & U
 \end{smallmatrix}\right]\in \Lambda$ we have that:
 
\begin{align*}
 & \Big[d_{\Ho_{\Dg(\T)} (A,A')}(\alpha)\coprod d_{\Ho_{\Dg(\U)} (B,B')}(\beta)\Big]_{\left[ \begin{smallmatrix}
T & 0\\
M & U
 \end{smallmatrix}\right]}=\\
 & =[d_{\Ho_{\Dg(\T)} (A,A')}(\alpha)]_{T}\coprod [d_{\Ho_{\Dg(\U)} (B,B')}(\beta)]_{U}
\end{align*}
By Equation \ref{dgNatmod}, we have that:
\begin{align*}
& [d_{\Ho_{\Dg(\T)} (A,A')}(\alpha)]_{T}:=d_{\mathrm{Hom}_{\mathrm{DgMod}(K)}(A(T),A'(T))}(\alpha_{T})=\\
& =d_{A'(T)}\circ \alpha_{T}-(-1)^{|\alpha_{T}|}\alpha_{T}\circ d_{A(T)}.
\end{align*}
Also, $[d_{\Ho_{\Dg(\U)} (B,B')}(\alpha)]_{U}:=d_{\mathrm{Hom}_{\mathrm{DgMod}(K)}(B(U),B'(U))}(\beta_{U})=d_{B'(U)}\circ \beta_{U}-(-1)^{|\beta_{U}|}\beta_{U}\circ d_{B(U)}$. Hence we conclude that
\begin{align*}
 \mathfrak{F}\Big(d_{\Delta_{1}}(\alpha,\beta)\Big)_{\left[ \begin{smallmatrix}
T & 0\\
M & U
 \end{smallmatrix}\right]}=\Big(d_{\Delta_{2}}(\mathfrak{F}(\alpha,\beta))\Big)_{\left[ \begin{smallmatrix}
T & 0\\
M & U
 \end{smallmatrix}\right]}.
\end{align*}
This proves that  the required diagram commutes, that is, $\mathfrak{F}$ commute with the differentials and hence
$\mathfrak{F}$ is a dg functor.
\end{proof}

\begin{lemma}\label{suprayectividad de funtor F}
Let $(A,f,B),(A',f',B')\in \Big( \mathrm{Mod}(\mathcal{T}),\mathbb{G}\mathrm{Mod}(\mathcal{U})\Big)$ be.
The map 
$$\mathfrak{F}:\mathrm{Hom}\Big((A,f,B),(A',f',B')\Big)\rightarrow \mathrm{Hom}_{\mathsf{Mod}(\Lambda)}(A\amalg _f B, A'\amalg _{f'} B')$$
is bijective. That is, $\mathfrak{F}$ is full and  faithful.
\end{lemma}
\begin{proof}
Firstly, we will see that $\mathfrak{F}$ is surjective.\\
Let $(A,f,B),(A',f',B')\in \Big( \mathrm{Mod}(\mathcal{T}),\mathbb{G}\mathrm{Mod}(\mathcal{U})\Big)$ and let $S:A\amalg_{f} B\longrightarrow A'\amalg_{f'} B'$
be a morphism in $\mathrm{Mod}(\Lambda)$ of degree $n$ whose components are:
$$S=\left\{S_{\left[\begin{smallmatrix}
T& 0 \\
M & U \\
  \end{smallmatrix} \right]}:A(T)\amalg B(U)\rightarrow A'(T)\amalg B'(U)\right\}_{\left[
\begin{smallmatrix}
T&0\\
M&U
\end{smallmatrix}
\right] \in \Lambda.}$$
That is, for each $\left[
\begin{smallmatrix}
T&0\\
M&U
\end{smallmatrix}
\right] \in \Lambda$ we have that:
$$S_{\left[\begin{smallmatrix}
T& 0 \\
M & U \\
  \end{smallmatrix} \right]}\in \mathrm{Hom}_{\mathrm{DgMod}(K)}^{n}\Big(A(T)\amalg B(U),A'(T)\amalg B'(U)\Big).$$
For $T\in \mathcal{T}$, consider the object $\left[\begin{smallmatrix}
T& 0 \\
M & 0 \\
\end{smallmatrix} \right]\in \Lambda$. Then we have the
morphism of degree $n$:
$$S_{\left[\begin{smallmatrix}
T& 0 \\
M & 0 \\
\end{smallmatrix} \right]}:A(T)\amalg 0\rightarrow A'(T)\amalg 0.$$
Therefore,  for $(x,0)\in (A(T))^{k}\amalg 0^{k}$  homogenous of degree $k$ we have that
$S_{\left[\begin{smallmatrix}
T& 0 \\
M & 0 \\
\end{smallmatrix} \right]}(x,0)=(x',0),$ for some $x'\in (A'(T))^{k+n}$.\\
Hence, we define $\alpha_{T}:A(T)\longrightarrow A'(T)$ as
$\alpha_{T}(x):=x'$ for $x\in A(T)$.
It is easy to see that $\alpha=\{\alpha_{T}:A(T)\longrightarrow A'(T)\}_{T\in \mathcal{T}}$ is a dg-natural transformation of degree $n$.\\
Now, for $U\in \mathcal{U}$ we consider the morphism $S_{\left[\begin{smallmatrix}
0& 0 \\
M & U \\
\end{smallmatrix} \right]}:0\amalg B(U)\rightarrow 0 \amalg B'(U).$
Then, we define $\beta_{U}:B(U)\longrightarrow B'(U)$ as follows:
for $y\in B(U)$ we have that $\beta_{U}(y):=y'\in B'(U)$ is such that $S_{\left[\begin{smallmatrix}
0& 0 \\
M & U \\
\end{smallmatrix} \right]}(0,y)=(0,y')$.
Similarly, we can prove that $\beta=\{\beta_{U}:B(U)\longrightarrow B'(U)\}$ is a
morphism of $\mathcal{U}$-modules of degree $n$.\\
Therefore, for $T\in \mathcal{T}$ and $U\in \mathcal{U}$ we define
$\alpha_{T}\amalg \beta_{U}:A(T)\amalg B(U)\longrightarrow A'(T)\amalg B'(U)$
as $(\alpha_{T}\amalg \beta_{U})(x,y):=(\alpha_{T}(x),\beta_{U}(y))$.
In this way, we obtain a family of morphisms in $\mathrm{DgMod}(K)$ as follows:
$$\alpha\amalg \beta:=\Big \{(\alpha_{T}\amalg \beta_{U}):A(T)\amalg B(U)\rightarrow A'(T)\amalg B'(U)\Big \}_{\left[
\begin{smallmatrix}
T&0\\
M&U
\end{smallmatrix}
\right] \in \Lambda.}$$
It is straightforward to see that $S=\alpha\amalg \beta$.\\
Now, let us check that $(\alpha,\beta)$ is a morphism from $(A,f,B)$ to $(A',f',B')$. 
We have to show that for $T\in \mathcal{T}$ the following diagram commutes in $\mathrm{DgMod}(K)$:
\[
\begin{diagram}
\node{A(T)} \arrow{e,t}{\alpha_{T}}\arrow{s,l}{f_{T}}\node{A'(T)}\arrow{s,r}{f'_{T}}\\
\node{\mathbb{G}(B)(T)} \arrow{e,b}{\mathbb{G}(\beta)_{T}}\node{\mathbb{G}(B')(T).}
\end{diagram}
\]
Then, for $x\in A(T)$ and homogenous element $m\in M^{p}(U,T)$ we have to show that
$$([\beta]_{U}\circ [f_{T}(x)]_{U})(m)=([f'_{T}(\alpha_{T}(x))]_{U})(m).$$

Consider an element $m \in M^{p}(U,T)$  , $0 \in \Ho_{\T}^{p}(T,T)$, $0 \in \Ho_{\U}^{p}(U,U)$. Then we have the morphism $\left[\begin{smallmatrix}
0 & 0\\
m & 0 \end{smallmatrix}\right]:\left[\begin{smallmatrix}
T & 0\\
M & U \end{smallmatrix}\right]\longrightarrow \left[\begin{smallmatrix}
T & 0\\
M & U\end{smallmatrix}\right]$ of degree $p$.\\

Since $S:A\amalg_{f}B\longrightarrow A'\amalg_{f'}B'$ is a dg-natural transformation in  $\mathrm{Mod}(\Lambda)$ of degree $n$, we have the following commutative diagram up to a sign $(-1)^{pn}$:

$$\xymatrix{A(T)\amalg B(U)\ar[rr]^{\alpha_{T}\amalg \beta_{U}}\ar[d]_{{\left[ \begin{smallmatrix}
0 & 0 \\
m & 0\\
\end{smallmatrix} \right]}}  & & A'(T)\amalg B'(U)\ar[d]^{{\left[ \begin{smallmatrix}
0 & 0 \\
m & 0\\
\end{smallmatrix} \right]}}\\
A(T)\amalg B(U)\ar[rr]^{\alpha_{T}\amalg \beta_{U}}
 & & A'(T)\amalg B'(U).}$$

Hence, on one hand, for each $(x,y) \in A(T)\amalg B(U)$ homogenous of degree $|(x,y)|=|x|=|y|$  we obtain that
 \begin{align*}
 \big( (\alpha_{T}\amalg \beta_{U})\circ \left[\begin{smallmatrix}
0 & 0\\
m & 0\end{smallmatrix}\right] \big)\begin{bmatrix}
x\\
y \end{bmatrix} = (\alpha_{T}\amalg \beta_{U}) \big( 0 , m \cdot x \big)
& = \big(\alpha_{T}(0), \beta_{U}(m \cdot x ) \big)\\
& = ( 0,\beta_{U}(m \cdot x ) ).
 \end{align*}
 
On the other hand, $\big( \left[\begin{smallmatrix}
0 & 0\\
m & 0 \end{smallmatrix}\right] \circ (\alpha_{T} \amalg \beta_{U})   \big) \begin{bmatrix}
x\\
y\end{bmatrix}  = \left[\begin{smallmatrix}
0 & 0\\
m & 0\end{smallmatrix}\right] \begin{bmatrix}
\alpha_{T}(x)\\
\beta_{U}(y)\end{bmatrix} = (0, m \cdot (\alpha_{T}(x)) )$. 
Consequently,  for $m \in M^{p}(U,T)$ we have that $\beta_{U}(m \cdot x )= (-1)^{pn} m \cdot \alpha_{T}(x).$

Now, we have the following $\beta_{U}(m \cdot x )= (-1)^{|x||m|} \beta_{U} \Big([f_{T}(x)]_{U}(m)\Big)$
and also we get that $m \cdot (\alpha_{T}(x)) = (-1)^{pn}(-1)^{|m|||x|} [f'_{T}(\alpha_{T}(x))]_{U}(m)$.
Then, from the previous equalities we have that:
$(-1)^{|x||m|} \beta_{U}\Big( [f_{T}(x)]_{U}(m)\Big) = (-1)^{|x||m|} [f'_{T}(\alpha_{T}(x))]_{U}(m).$
Therefore, for  $m \in M^{p}(U,T)$  we obtain that:
$$\beta_{U} \Big([f_{T}(x)]_{U}(m)\Big)=  [f'_{T}(\alpha_{T}(x))]_{U}(m).$$
Thus, $f'_{T}(\alpha_{T}(x))= \beta \circ f_{T}(x)= \Ho_{\Dg(\U)}(M_{T}, \beta)(f_{T}(x))= (\mathbb{G}(\beta)_{T})(f_{T}(x))$ for each  homogenous element $x \in A(T)$. Thus, $f'_{T}\circ \alpha_{T}=\mathbb{G}(\beta)_{T}\circ f_{T}$. Therefore,
$(\alpha, \beta): (A,f,B)\longrightarrow (A',f',B')$ is a morphism in $\big( \Dg(\T),\mathbb{G} \big(\Dg(\U)\big)\big)$ and  $\mathfrak{F}(\alpha, \beta)= S$.  We conclude that $\mathfrak{F}$ is full.\\
Now, let suppose that $(\alpha',\beta'),(\alpha,\beta):(A,f,B)\longrightarrow (A',f',B')$  are morphisms
such that and $\mathfrak{F}(\alpha,\beta)=\mathfrak{F}(\alpha',\beta')$. Then for each $T\in \mathcal{T}$, $U\in \mathcal{U}$
we have that $\alpha_{T}\amalg \beta_{U}=\alpha'_{T}\amalg \beta'_{U}$. This implies, that
$\alpha_{T}=\alpha'_{T}$ and $\beta_{U}=\beta'_{U}$ and therefore $\alpha=\alpha'$ and
$\beta=\beta'$. Proving that $\mathfrak{F}$ is injective and therefore $\mathfrak{F}$ is full and faithful.
\end{proof}

 We can define  a dg-functor $I_{1}: \T \longrightarrow \Lambda$  as follows $I_{1}(T):=\left[\begin{smallmatrix}
T & 0\\
M & 0 \end{smallmatrix}\right] $ and for a homogenous morphism  $t:T \longrightarrow T'$ of degree $|t|$ in $\T$ we set
$$I_{1}(t)=\left[\begin{smallmatrix}
t & 0 \\
0 & 0
\end{smallmatrix} \right]: \left[\begin{smallmatrix}
T & 0 \\
M & 0
\end{smallmatrix} \right] \longrightarrow \left[\begin{smallmatrix}
T' & 0 \\
M & 0
\end{smallmatrix} \right].
$$
Then we have that $I_{1}$ is a graded functor. For simplicity let us denote $\Delta_{1}:=\mathrm{Hom}_{\Lambda}\Big(\left[\begin{smallmatrix}
T & 0 \\
M & 0
\end{smallmatrix} \right] ,\left[\begin{smallmatrix}
T' & 0 \\
M & 0
\end{smallmatrix} \right] \Big)$. It is easy to see that the following diagram commutes

$$\xymatrix{\mathrm{Hom}_{\mathcal{T}}(T,T')\ar[rr]^{I_{1}}\ar[d]_{d_{\mathcal{T}(T,T')}} & &{ \mathrm{Hom}_{\Lambda}\Big(\left[\begin{smallmatrix}
T & 0 \\
M & 0
\end{smallmatrix} \right] ,\left[\begin{smallmatrix}
T' & 0 \\
M & 0
\end{smallmatrix} \right] \Big)}\ar[d]^{d_{\Delta_{1}}}\\
\mathrm{Hom}_{\mathcal{T}}(T,T')\ar[rr]^{I_{1}} & & {\mathrm{Hom}_{\Lambda}\Big(\left[\begin{smallmatrix}
T & 0 \\
M & 0
\end{smallmatrix} \right] ,\left[\begin{smallmatrix}
T' & 0 \\
M & 0
\end{smallmatrix} \right] \Big).}}$$ Hence,  $I_{1}$ is a dg-functor. In the same way, we define a dg-functor $I_{2}: \U \longrightarrow \Lambda$. Then we have the induced  dg-functors 
$$\mathbb{I}_{1}: \Dg(\Lambda) \longrightarrow \Dg(\T)$$
$$\mathbb{I}_{2}: \Dg(\Lambda) \longrightarrow \Dg(\U).$$
For $C$ a dg $\Lambda$-module, we denote by $C_{1}:= \mathbb{I}_{1}(C)= C \circ I_{1}: \T \longrightarrow \Dg(K)$ where $ C_{1}$ is given as follows:
$$C_{1}(T)=C\big(  \left[\begin{smallmatrix}
T & 0\\
M & 0\end{smallmatrix}\right] \big)$$
$$C_{1}(t)=C \big(  \left[\begin{smallmatrix}
t & 0\\
0 & 0\end{smallmatrix}\right] \big)$$
and $ C_{2}= \mathbb{I}_{2}(C)= C \circ I_{2}: \U \longrightarrow \Dg(K)$  is defined analogously.

\begin{lemma}\label{lema1 equ c1 a Gc2}
Let $C$ be a  dg $\Lambda$-module. Then, there exists a
morphism of dg $\mathcal{T}$-modules of degree zero
$$f:C_{1}\longrightarrow \mathbb{G}(C_{2}),$$
such that $D_{\mathrm{Hom}_{\mathrm{DgMod}(\T)}(C_{1},\mathbb{G}(C_{2}))}(f)=0$.
\end{lemma}
\begin{proof}
Let $C$ be a dg $\Lambda$-module and consider
$T\in\mathcal T$ and $U\in\mathcal U$. For all $m\in M(U,T)$ homogenous of degree $|m|$
we have  $\overline{m}:=\left[
\begin{smallmatrix}
0& 0 \\
m & 0 \\
\end{smallmatrix}
\right]:
\left[
\begin{smallmatrix}
T& 0 \\
M & 0 \\
\end{smallmatrix}
\right]\rightarrow
\left[
\begin{smallmatrix}
0& 0 \\
M & U \\
\end{smallmatrix}
\right]$ of degree $|m|$ in $\Lambda$.
We note that $\left[
\begin{smallmatrix}
T& 0 \\
M & 0 \\
\end{smallmatrix}
\right]=I_{1}(T)$ and $\left[
\begin{smallmatrix}
0& 0 \\
M & U \\
\end{smallmatrix}
\right]=I_{2}(U)$
Applying the  dg $\Lambda$-module $C$  to $\overline{m}$ yields a morphism in $\mathrm{DgMod}(K)$ of degree $|m|$:
$$C(\overline{m}):C_{1}(T)=C(I_{1}(T))\longrightarrow C_{2}(U)=C(I_{2}(U)).$$
We assert that the morphisms $C(\overline{m})$ induces a morphism of degree zero in $\mathrm{DgMod}(K)$:
$$f_{T}:C_{1}(T)\longrightarrow \mathbb{G}(C_{2})(T)=\mathrm{Hom}_{\mathrm{Mod}(\mathcal{U})}(M_{T},C_{2}).$$
Indeed, for  a homegenous element $x\in C_{1}(T) $ we define $f_{T}(x):M_{T}\longrightarrow C_{2}$ with 
$$f_{T}(x)=\left\lbrace [f_{T}(x)]_{U}:M_{T}(U)\longrightarrow C_{2}(U) \right\rbrace _{U\in \mathcal{U}},$$ 
where $[f_{T}(x)]_{U}(m):=(-1)^{|x||m|}C(\overline{m})(x)$ for all  homogenous element $m\in M_{T}(U)$. Let us check that
$f_{T}(x):M_{T}\longrightarrow C_{2}$ is a dg natural transformation of degree $|x|$.
Indeed, consider $U\in \mathcal{U}$. Hence, for a homogenous element $m\in M(U,T)$ of degree $|m|$ we have that
$[f_{T}(x)]_{U}(m)=(-1)^{|x||m|}C(\overline{m})(x)\in C_{2}(U)^{|x|+|m|}$ since $C(\overline{m})$ is of degree $|x|$. This proves that $[f_{T}(x)]_{U}$ is a graded map of degree $|x|$ and hence $f_{T}(x):M_{T}\longrightarrow C_{2}$ is a graded map of degree $|x|$.
Now, let $u\in\mathrm{Hom}_{\mathcal{U}}(U,U')$ a homogenous element of degree $|u|$, is is easy to show that the following diagram commutes up to the sign $(-1)^{|x||u|}$

\[
\begin{diagram}
\node{M_{T}(U)}\arrow{e,t}{[f_{T}(x)]_{U}}\arrow{s,l}{M_{T}(u)}
 \node{C_{2}(U)}\arrow{s,r}{C_{2}(u)}\\
\node{M_{T}(U')}\arrow{e,b}{[f_{T}(x)]_{U'}}
 \node{C_{2}(U').}
\end{diagram}
\]
 Hence $f_{T}(x) $ is a morphism of dg $ \mathcal{U}$-modules of degree $|x|$. This proves that $f_{T}:C_{1}(T)\longrightarrow \mathbb{G}(C_{2})(T)=\mathrm{Hom}_{\mathrm{Mod}(\mathcal{U})}(M_{T},C_{2})$
is a graded map of degree $0$.\\
Now, let $ t\in \mathrm{Hom}_{\mathcal{T}}(T,T')$ a homogenous element, we assert that the following diagram commutes
\[
\begin{diagram}
\node{C_{1}(T)}\arrow{e,t}{f_{T}}\arrow{s,l}{C_{1}(t)}
 \node{\mathbb{G}(C_{2})(T)}\arrow{s,r}{\mathbb{G}(C_{2})(t)}\\
\node{C_{1}(T')}\arrow{e,b}{f_{T'}}
 \node{\mathbb{G}(C_{2})(T').}
\end{diagram}
\]
Indeed, let $x\in C_{1}(T)$ an homogenous element.  Then
\begin{eqnarray*}
\Big(\mathbb{G}(C_{2})(t)\circ f_{T}\Big)(x) = \Big(\mathrm{Hom}_{\mathrm{DgMod}(\mathcal{U})}(\bar{t},C_{2}) \circ f_{T}\Big)(x)= (-1)^{|t||x|}f_{T}(x)\circ \bar{t}.
\end{eqnarray*}
For $U\in \mathcal{U}$ and a homogenous element $m'\in M(U,T')$ we have that
$$(-1)^{|t||x|}\Big([f_{T}(x)]_{U}\circ[\bar{t}]_{U}\Big)(m')=(-1)^{|m'||t|+|m'||x|}C\Big(\left[
\begin{smallmatrix}
0&0\\
m'\bullet t &0
\end{smallmatrix}
\right]\Big)(x).$$
On the other hand, $(f_{T'}\circ C_{1}(t))(x)=f_{T'}(C_{1}(t)(x))
=f_{T'}\Big(C\Big(\left[
\begin{smallmatrix}
t&0\\
0 &0
\end{smallmatrix}
\right]\Big)(x)\Big )$. Then for $U\in \mathcal{U}$ and $m'\in M(U,T')$ we have that:
\begin{eqnarray*}
\Big[f_{T'}\Big(C\Big(\left[
\begin{smallmatrix}
t&0\\
0 &0
\end{smallmatrix}
\right]\Big)(x)\Big )\Big]_{U}(m')  =  (-1)^{e\cdot |m'|}C\Big(\left[
\begin{smallmatrix}
0&0\\
m'\bullet t &0
\end{smallmatrix}
\right]\Big)(x),
\end{eqnarray*}
where $e=|C\Big(\left[
\begin{smallmatrix}
t&0\\
0 &0
\end{smallmatrix}
\right]\Big)(x)|=|t|+|x|$ (this happens since $C$ is a graded morphism of degree zero and hence $C\Big(\left[
\begin{smallmatrix}
t&0\\
0 &0
\end{smallmatrix}
\right]\Big)$ has degree $|t|$). Hence, $\Big(\mathbb{G}(C_{2})(t)\circ f_{T}\Big)(x)=\Big(f_{T'}\circ C_{1}(t)\Big)(x)$ and the required diagram commutes.\\
Now, let us see that $D(f)=0$, where
$$D:\mathrm{Hom}_{\mathrm{DgMod}(\T)}(C_{1},\mathbb{G}(C_{2}))\longrightarrow \mathrm{Hom}_{\mathrm{DgMod}(\T)}(C_{1},\mathbb{G}(C_{2}))$$ is defined for each $T\in \mathcal{T}$ as follows:  $(D(f))_{T}:=d(f_{T})$ (see Equation \ref{dgNatmod}). Recall that $d:\mathrm{Hom}_{\mathrm{DgMod}(K)}(C_{1}(T),\mathbb{G}(C_{2})(T))\longrightarrow \mathrm{Hom}_{\mathrm{DgMod}(K)}(C_{1}(T),\mathbb{G}(C_{2})(T))$ is defined as
$d(\alpha)=d_{\mathbb{G}(C_{2})(T)}\circ \alpha-(-1)^{|\alpha|}\alpha\circ d_{C_{1}(T)}$ for every homogenous element $\alpha \in \mathrm{Hom}_{\mathrm{DgMod}(K)}(C_{1}(T),\mathbb{G}(C_{2})(T))$ (see Equation \ref{Homdgestructure}). Then in order to see that $D(f)=0$, it is enough to see that for each $T\in \mathcal{T}$ the following diagram commutes
$$(\ast):\xymatrix{C_{1}(T)\ar[r]^{f_{T}}\ar[d]_{d_{C_{1}(T)}} & G(C_{2})(T)\ar[d]^{d_{\mathbb{G}(C_{2})(T)}}\\
C_{1}(T)\ar[r]^{f_{T}} & G(C_{2})(T).}$$
Then, in order to show that the diagram $(\ast)$ commutes, we have to show that for each homogenous element $x\in C_{1}(T)=C\Big(\left[
\begin{smallmatrix}
T&0\\
M &0
\end{smallmatrix}
\right]\Big)$ and $U\in \mathcal{U}$ we have that
$$\Big[d_{\mathbb{G}(C_{2})(T)}\Big(f_{T}(x)\Big)\Big]_{U}=\Big[f_{T}\Big(d_{C_{1}(T)}(x)\Big)\Big]_{U}.$$
Recall that, by Equations \ref{Homdgestructure} and \ref{dgNatmod}, we have that
$$d_{\mathbb{G}(C_{2})(T)}:\mathrm{Hom}_{\mathrm{DgMod}(\mathcal{U})}(M_{T},C_{2})\longrightarrow \mathrm{Hom}_{\mathrm{DgMod}(\mathcal{U})}(M_{T},C_{2})$$
is defined for $\gamma\in \mathrm{Hom}_{\mathrm{DgMod}(\mathcal{U})}(M_{T},C_{2})$ as follows:
$(d_{\mathbb{G}(C_{2})(T)}(\gamma))_{U}:=\delta (\gamma_{U})=d_{C\Big(\left[
\begin{smallmatrix}
0&0\\
M &U
\end{smallmatrix}
\right]\Big)}\circ \gamma_{U}-(-1)^{|\gamma_{U}|}\gamma_{U}\circ d_{M(U,T)}$
for each $U\in \mathcal{U}$, where $\delta$ is the differential of de dg K-module $\mathrm{Hom}_{\mathrm{DgMod}(K)}\Big(M(U,T),C\Big(\left[
\begin{smallmatrix}
0&0\\
M &U
\end{smallmatrix}
\right]\Big)\Big)$.\\
On one hand, by the above discussion we obtain that
\begin{align*}
& \Big[d_{\mathbb{G}(C_{2})(T)}\Big(f_{T}(x)\Big)\Big]_{U} =d_{C\Big(\left[
\begin{smallmatrix}
0&0\\
M &U
\end{smallmatrix}
\right]\Big)}\circ [f_{T}(x)]_{U}-(-1)^{|[f_{T}(x)]_{U}|}[f_{T}(x)]_{U}\circ d_{M(U,T)}.
\end{align*}
On the other hand, we get that
\begin{align*}
\Big[f_{T}\Big(d_{C_{1}(T)}(x)\Big)\Big]_{U}=\left [f_{T}\Big(d_{C\Big(\left[
\begin{smallmatrix}
T&0\\
M & 0
\end{smallmatrix}
\right]\Big)}(x)\Big)\right]_{U}.
\end{align*}
In order to show that the two previous equalities coincide, we first recall that $C$ is a dg $\Lambda$-module and hence it implies that the following diagram commutes
$$\xymatrix{\mathrm{Hom}_{\Lambda}{\Big(\left[
\begin{smallmatrix}
T& 0 \\
M & 0 \\
\end{smallmatrix}
\right],
\left[
\begin{smallmatrix}
0& 0 \\
M & U \\
\end{smallmatrix}
\right]\Big)}\ar[rr]^{C}\ar[d]^{\Delta_{1}}  & &\mathrm{Hom}_{\mathrm{DgMod}(K)}{\Big(C\Big(\left[
\begin{smallmatrix}
T& 0 \\
M & 0 \\
\end{smallmatrix}
\right]\Big),
C\Big(\left[
\begin{smallmatrix}
0& 0 \\
M & U \\
\end{smallmatrix}
\right]\Big)\Big)}\ar[d]^{\Delta_{2}}\\
\mathrm{Hom}_{\Lambda}{\Big(\left[
\begin{smallmatrix}
T& 0 \\
M & 0 \\
\end{smallmatrix}
\right],
\left[
\begin{smallmatrix}
0& 0 \\
M & U \\
\end{smallmatrix}
\right]\Big)}\ar[rr]^{C}  & &\mathrm{Hom}_{\mathrm{DgMod}(K)}{\Big(C\Big(\left[
\begin{smallmatrix}
T& 0 \\
M & 0 \\
\end{smallmatrix}
\right]\Big),
C\Big(\left[
\begin{smallmatrix}
0& 0 \\
M & U \\
\end{smallmatrix}
\right]\Big)\Big),}}$$
where $\Delta_{1}$ and $\Delta_{2}$ are the corresponding differentials. Then, for $\left[
\begin{smallmatrix}
0& 0 \\
m & 0 \\
\end{smallmatrix}
\right]:
\left[
\begin{smallmatrix}
T& 0 \\
M & 0 \\
\end{smallmatrix}
\right]\rightarrow
\left[
\begin{smallmatrix}
0& 0 \\
M & U \\
\end{smallmatrix}
\right]$ we have that

\begin{align*}
C\Big(\left[
\begin{smallmatrix}
0& 0 \\
d_{M(U,T)}(m) & 0 \\
\end{smallmatrix}
\right]\Big) & =d_{C\Big(\left[
\begin{smallmatrix}
0& 0 \\
M & U \\
\end{smallmatrix}
\right]\Big)}\circ C\Big(\left[
\begin{smallmatrix}
0& 0 \\
m & 0 \\
\end{smallmatrix}
\right]\Big)-(-1)^{|m|}C\Big(\left[
\begin{smallmatrix}
0& 0 \\
m & 0 \\
\end{smallmatrix}
\right]\Big)\circ d_{C\Big(\left[
\begin{smallmatrix}
T& 0 \\
M & 0\\
\end{smallmatrix}
\right]\Big)}.
\end{align*}
By definition of $[f_{T}(x)]_{U}\Big(d_{M(U,T)}(m)\Big)$, it follows that:

\begin{align*}
& [f_{T}(x)]_{U}\Big(d_{M(U,T)}(m)\Big)=\\
 &=(-1)^{|x|\cdot |d_{M(U,T)}(m)|}
C\Big(\left[
\begin{smallmatrix}
0& 0 \\
d_{M(U,T)}(m) & 0 \\
\end{smallmatrix}
\right]\Big)(x)\\
& =(-1)^{|x|\cdot (|m|+1)}\Big(d_{C\Big(\left[
\begin{smallmatrix}
0& 0 \\
M & U \\
\end{smallmatrix}
\right]\Big)}\circ C\Big(\left[
\begin{smallmatrix}
0& 0 \\
m & 0 \\
\end{smallmatrix}
\right]\Big)-(-1)^{|m|}C\Big(\left[
\begin{smallmatrix}
0& 0 \\
m & 0 \\
\end{smallmatrix}
\right]\Big)\circ d_{C\Big(\left[
\begin{smallmatrix}
T& 0 \\
M & 0\\
\end{smallmatrix}
\right]\Big)}\Big)(x).
\end{align*}
Therefore we obtain that:

\begin{align*}
& (-1)^{|x|}[f_{T}(x)]_{U}\Big(d_{M(U,T)}(m)\Big)=\\
& =\!(-1)^{|x||m|}\Big(\!d_{C\Big(\left[
\begin{smallmatrix}
0& 0 \\
M & U \\
\end{smallmatrix}
\right]\Big)}\!\!\circ\! C\Big(\left[
\begin{smallmatrix}
0& 0 \\
m & 0 \\
\end{smallmatrix}
\right]\Big)\!\Big)(x)\!-\!(-1)^{|m|+|x||m|}C\Big(\!\left[
\begin{smallmatrix}
0& 0 \\
m & 0 \\
\end{smallmatrix}
\right]\Big)\!\!\circ d_{C\Big(\left[
\begin{smallmatrix}
T& 0 \\
M & 0\\
\end{smallmatrix}
\right]\Big)}\!\Big)(x)\\
&= d_{C\Big(\!\left[
\begin{smallmatrix}
0& 0 \\
M & U \\
\end{smallmatrix}
\right]\!\Big)}\Big(\![f_{T}(x)]_{U}(m)\Big)\!-\!(-1)^{|m|+|x||m|}(-1)^{|m|(|x|+1)}\!\Big [f_{T}\Big(d_{C\Big(\!\left[
\begin{smallmatrix}
T&0\\
M & 0
\end{smallmatrix}
\right]\!\Big)}(x)\Big)\!\Big]_{U}\!(m)\\
&= d_{C\Big(\left[
\begin{smallmatrix}
0& 0 \\
M & U \\
\end{smallmatrix}
\right]\Big)}\Big([f_{T}(x)]_{U}(m)\Big)-\Big [f_{T}\Big(d_{C\Big(\left[
\begin{smallmatrix}
T&0\\
M & 0
\end{smallmatrix}
\right]\Big)}(x)\Big)\Big]_{U}(m).
\end{align*}
Hence, we conclude that
\begin{align*}
\Big [f_{T}\Big(d_{C\Big(\left[
\begin{smallmatrix}
T&0\\
M & 0
\end{smallmatrix}
\right]\Big)}(x)\Big)\Big]_{U}& = d_{C\Big(\left[
\begin{smallmatrix}
0& 0 \\
M & U \\
\end{smallmatrix}
\right]\Big)}\circ [f_{T}(x)]_{U}-(-1)^{|[f_{T}(x)]_{U}|}[f_{T}(x)]_{U}\circ d_{M(U,T)}.
\end{align*}
Proving that
$\Big[d_{\mathbb{G}(C_{2})(T)}\Big(f_{T}(x)\Big)\Big]_{U}= \Big[f_{T}\Big(d_{C_{1}(T)}(x)\Big)\Big]_{U};$
and hence the diagram $(\ast)$ commutes, this shows that $D(f)=0$. 
Hence $\left(C_{1},f,C_{2} \right)\in \Big(\mathrm{Mod}(\mathcal{T}),\mathbb{G}\big(\mathrm{Mod}(\mathcal{U})\big) \Big)$.
\end{proof}

\begin{lemma}\label{C iso a C1 C2}
Let $C$ be a $\Lambda$ dg k-module. Then
$C_1\underset{f}\amalg C_2\cong C$.
\end{lemma}
\begin{proof}
Let $T\in\mathcal T$ and $U\in\mathcal U$.  We have sequences of morphism of degree zero in $\Lambda$.
$$\left[ {\begin{array}{cc}
T& 0 \\
M& 0 \\
\end{array} } \right]\xrightarrow{\lambda_T:=\left[ {\begin{smallmatrix}
1_T& 0 \\
0& 0 \\
\end{smallmatrix} } \right]} \left[ {\begin{array}{cc}
T& 0 \\
M& U \\
\end{array} } \right]\xrightarrow{\rho_T:=\left[ {\begin{smallmatrix}
1_T& 0 \\
0& 0 \\
\end{smallmatrix} } \right]} \left[ {\begin{array}{cc}
T& 0 \\
M& 0 \\
\end{array} } \right]$$

$$\left[ {\begin{array}{cc}
0& 0 \\
M& U \\
\end{array} } \right]\xrightarrow{\lambda_U:=\left[ {\begin{smallmatrix}
0& 0 \\
0& 1_U\\
\end{smallmatrix} } \right]} \left[ {\begin{array}{cc}
T& 0 \\
M& U \\
\end{array} } \right]\xrightarrow{\rho_U:=\left[ {\begin{smallmatrix}
0& 0 \\
0& 1_U \\
\end{smallmatrix} } \right]} \left[ {\begin{array}{cc}
0& 0 \\
M& U \\
\end{array} } \right]$$

Since $C$ is a $\Lambda$ dg module we get the following maps of degree zero in $\mathrm{DgMod}(K)$:

$$C(\lambda_{T}):=C\Big(\left[
\begin{smallmatrix}
1_{T}&0\\
0 &0
\end{smallmatrix}
\right]\Big): C\Big(\left[
\begin{smallmatrix}
T&0\\
M & 0
\end{smallmatrix}
\right]\Big)=C_{1}(T)\longrightarrow
C\Big(\left[
\begin{smallmatrix}
T &0\\
M & U
\end{smallmatrix}
\right]\Big)$$

$$C(\lambda_{U}):=C\Big(\left[
\begin{smallmatrix}
0 &0\\
0 & 1_{U}
\end{smallmatrix}
\right]\Big): C\Big(\left[
\begin{smallmatrix}
0&0\\
M & U
\end{smallmatrix}
\right]\Big)=C_{2}(U)\longrightarrow
C\Big(\left[
\begin{smallmatrix}
T &0\\
M & U
\end{smallmatrix}
\right]\Big).$$
This produces a map (by the universal poperty of the coproduc)
$$\phi_{T,U}:C_{1}(T)\amalg C_{2}(U)=\Big(C_{1}\amalg_{f}C_{2}\Big)\Big(\left[
\begin{smallmatrix}
T &0\\
M & U
\end{smallmatrix}
\right]\Big)\longrightarrow
C\Big(\left[
\begin{smallmatrix}
T &0\\
M & U
\end{smallmatrix}
\right]\Big)$$
defined as $\phi_{T,U}(x,y)=C(\lambda_{T})(x)+C(\lambda_{U})(y)$ $\forall (x,y)\in C_{1}(T)\amalg C_{2}(U)$.\\
It is easy to see that $\phi=\!\left\{\!\!
\phi_{{}_{\left[
\begin{smallmatrix}
T& 0 \\
M& U\\
\end{smallmatrix}
\right] }}\!\!:=\!\phi_{T,U}\!\!:\!\! \Big(C_{1}\amalg_{f} C_{2}\Big) (\left [\begin{smallmatrix}
T& 0 \\
M& U\\
\end{smallmatrix} \right])\rightarrow  C\Big(\left [\begin{smallmatrix}
T& 0 \\
M& U\\
\end{smallmatrix} \right]\Big)\!\!\right\}_{\!\!\left[\begin{smallmatrix}
T& 0 \\
M& U\\
\end{smallmatrix} \right ]\in \Lambda}$
 is a dg-natural transformation of degree zero $\phi:C_{1}\amalg_{f} C_{2}\longrightarrow C$.
 Finally, it is easy to see that $\phi$ is invertible.
 \end{proof}

\begin{theorem}\label{equivalencia Lambda y coma}
The differential graded functor
$$\mathfrak{F}: \Big( \Dg(\T), \big(\mathbb{G}\Dg(\U)\big) \Big) \longrightarrow \Dg(\Lambda) $$
is a dg-equivalence of dg-categories.
\end{theorem}
\begin{proof}
It follows  from Lemmas  \ref{suprayectividad de funtor F}, \ref{lema1 equ c1 a Gc2} and  \ref{C iso a C1 C2}.
\end{proof}

Let $A$ be a dg-K-algebra. Then  there exists a dg-$K$-category $\mathcal{A}$ with only one object, say $\ast$, and 
$\mathrm{Hom}_\mathcal{A}(\ast,\ast)=A$. Then, the  evaluation  functor
$\mathrm{ev}_{\ast}: \mathrm{DgMod}(\mathcal{A})\rightarrow \mathrm{DgMod}(A)$ given as $\mathrm{ev}_{\ast}(F)=F(\ast)$, induces an equivalence of categories.\\
Let $U$ and $T$ be two dg-$K$-algebras and consider their respective dg-categories $\mathcal{U}$ and $\mathcal{T}$. Denote by 
$\ast$ and $\star$ the unique objects in  $\mathcal{U}$ and $\mathcal{T}$ respectively, that is, 
$\mathrm{Hom}_{\mathcal{U}}(\ast,\ast)=U$ and $\mathrm{Hom}_{\mathcal{T}}(\star,\star)=T$. Thus, by the above mentioned we have equivalences
of categories $\mathrm{DgMod}(\mathcal{U})\cong\mathrm{DgMod}(U)$ and $\mathrm{DgMod}(\mathcal{T})\cong\mathrm{DgMod}(T)$. In addition, 
the  dg-category tensor product $\mathcal{U}\otimes_{K}\mathcal{T}^{op}$ has only one object, the pair
$(\ast,\star)$, and $\mathrm{Hom}_{\mathcal{U}\otimes_{K}\mathcal{T}^{op}}((\ast,\star),(\ast,\star))=U\otimes_{K} T^{op}$. \\
In this way, we 
have the equivalence of categories $\mathrm{DgMod}(\mathcal{U}\otimes_{K}\mathcal{T}^{op})\cong\mathrm{DgMod}(U\otimes_{K} T^{op})$, and 
given 
a dg-$U$-$T$-bimodule  $N$, there exists a dg-functor 
$M:\mathcal U\otimes_{K} \mathcal T^{op}\rightarrow \mathrm{DgMod}(K)$  and an isomorphism of
 $U$-$T$-modules  $N\cong  M((\mathfrak{a},\star))$ such that,  after identifying $N$ with $ M((\ast,\star))$, 
$$u m=M(u\otimes 1_{\star})(m) \text{ and } m t = (-1)^{|m|\,|t|}M(1_{\ast}\otimes t^{op})(m),$$
 for all $u\in U$, $t\in T$ and  $m\in M((\ast,\star))$ homogenous elements.\\
 In this way, the dg-$K$-algebra $\Lambda=\left[
\begin{smallmatrix}
T& 0 \\ N & U
\end{smallmatrix}\right]$ produces a dg-$K$ category $\bf{\Lambda}=\left[
\begin{smallmatrix}
\mathcal{T} & 0 \\ M & \mathcal{U}
\end{smallmatrix}\right]$ with just one object $\left[
\begin{smallmatrix}
\star & 0 \\ M & \ast
\end{smallmatrix}\right]$ such that $\mathrm{Hom}_{\bf{\Lambda}}(\left[
\begin{smallmatrix}
\star & 0 \\ M & \ast
\end{smallmatrix}\right], \left[
\begin{smallmatrix}
\star & 0 \\ M & \ast
\end{smallmatrix}\right])=\Lambda$ and  the evaluation functor
$$e_{\left[
\begin{smallmatrix}
\star & 0 \\ M & \ast
\end{smallmatrix}\right]}:\mathrm{DgMod}({\bf{\Lambda}})\longrightarrow \mathrm{DgMod}(\Lambda)$$
is an equivalence of dg-$K$-categories. Now, consider the functor $\mathbb{G}:\mathrm{DgMod}(\mathcal{U})\longrightarrow \mathrm{DgMod}(\mathcal{T})$ constructed in Definition \ref{definiG}. Then, the following diagram commutes up isomorphism
 
$$\xymatrix{\mathrm{DgMod}(\mathcal {U})\ar[rr]^{\mathbb{G}}\ar[d]_{\mathrm{ev}_{\ast}} & &  \mathrm{DgMod}(\mathcal {T})\ar[d]^{{\mathrm{ev}_{\star}}}\\
\mathrm{DgMod}(U)\ar[rr]^{\mathrm{Hom}_{U}(N,-)} & & \mathrm{DgMod}(T).}$$

Indeed, let $B\in\mathrm{DgMod}(\mathcal{U})$ be. One  one hand 
$\mathrm{ev}_{\star}(\mathbb{G}(B))=\mathrm{Hom}_{\mathrm{DgMod}(\mathcal{U})}(M_{\star},B)$ and on the other hand we have 
$\mathrm{Hom}_{U}(N, \mathrm{ev}_{\ast}(B))=\mathrm{Hom}_{U}(N, B({\ast}))$.\\
Finally, the assertion follows from the isomorphism
$\mathrm{ev}_{\ast}:\mathrm{Hom}_{\mathrm{DgMod}(\mathcal{U})}(M_{\star},B)\rightarrow \mathrm{Hom}_U(M(\ast,\star), B(\ast))$, since $M(\ast,\star)\cong N$. Then, we can identify $\mathbb{G}$ with $\mathrm{Hom}_{U}(N,-)$.\\
Now, by Definition \ref{commadgcat} and Proposition \ref{dgcommaprop}, in this case $(\mathrm{DgMod}(T), \mathbb{G} \mathrm{DgMod}(U))$ is the dg-category whose objects are  the triples  $(A,f,B)$ with $A \in \mathrm{DgMod}(T)$, $B \in \mathrm{DgMod}(U)$ and $f:A \longrightarrow \mathrm{Hom}_{U}(N, B)$ a morphism of differential graded  $T$-modules of degree zero and such that $d_{\mathrm{Hom}_{T}(A,\mathrm{Hom}_{U}(N, B))}(f)=0$ and  given two objects  $(A,f,B)$ and  $(A',f',B')$ in  $(\mathrm{DgMod}(T), \mathbb{G} \mathrm{DgMod}(U))$ a morphism between this two objects  is a pair  of  homogenous morphisms  $(\alpha, \beta)$ with  $\alpha: A \longrightarrow A'$ a morphism of dg $T$-modules and  $\beta: B \longrightarrow  B$ a morphism  of dg $U$-modules  such $k=|\alpha|=|\beta|$ and such that  the following  diagram commutes 
 $$\xymatrix{ A \ar[rr]^{\alpha} \ar[d]_{f}& & A'  \ar[d]^{f'} \\
 \mathrm{Hom}_{U}(N, B) \ar[rr]_{\mathrm{Hom}_{U}(N, \beta)}& &  \mathrm{Hom}_{U}(N, B').}$$
 
Then we have the following result wich is a generalization of a well-known  result for artin algebras (see \cite[Proposition 2.2]{AusBook}), but now for dg-$K$-algebras.\\ 

\begin{corollary}
Let $U$ and $T$ two dg-K-algebras and $N$ a dg-$U$-$T$-bimodule. Then  $\Lambda=\left[
\begin{smallmatrix}
T& 0 \\ N & U
\end{smallmatrix}\right]$ is a dg-$K$-algebra and there exists an equivalence of dg-categories
$$(\mathrm{DgMod}(T), \mathbb{G} \mathrm{DgMod}(U))\cong \mathrm{DgMod}(\Lambda)$$
where $\mathbb{G}=\mathrm{Hom}_{U}(N, -)$. 
\end{corollary}


\footnotesize

\vskip3mm \noindent Martha Lizbeth Shaid Sandoval Miranda:\\ Departamento de Matem\'aticas, Universidad Aut\'onoma Metropolitana Unidad Iztapalapa\\
Av. San Rafael Atlixco 186, Col. Vicentina Iztapalapa 09340, M\'exico, Ciudad de M\'exico.\\ {\tt marlisha@xanum.uam.mx, marlisha@ciencias.unam.mx}

\vskip3mm \noindent Valente Santiago Vargas:\\ Departamento de Matem\'aticas, Facultad de Ciencias, Universidad Nacional Aut\'onoma de M\'exico\\
Circuito Exterior, Ciudad Universitaria,
C.P. 04510, Ciudad de M\'exico, MEXICO.\\ {\tt valente.santiago@ciencias.unam.mx}

\vskip3mm \noindent  Edgar Omar Velasco P\'aez:\\ Departamento de Matem\'aticas, Facultad de Ciencias, Universidad Nacional Aut\'onoma de M\'exico\\
Circuito Exterior, Ciudad Universitaria,
C.P. 04510, Ciudad de M\'exico, MEXICO.\\ {\tt edgar-bkz13@ciencias.unam.mx }

\end{document}